\documentclass[12pt]{article}

\usepackage{amsfonts, amssymb, epsfig,subfigure}
\usepackage{mathrsfs}
\usepackage{stmaryrd}
\usepackage{graphicx,amsmath}
\usepackage{latexsym}
\usepackage{verbatim}
\usepackage{color}

\newtheorem{thm}{Theorem}[section]

\newtheorem{assp}{Assumption}
\newtheorem{lem}[thm]{Lemma}

\newtheorem{cor}[thm]{Corollary}

\newtheorem{expl}[thm]{Example}

\makeatletter
\@addtoreset{equation}{section}

\newcommand{\beq}[1]{\begin{equation} \label{#1}}
\newcommand{\eeq}{\end{equation}}
\newcommand{\bed}{\begin{displaymath}}
\newcommand{\eed}{\end{displaymath}}
\newcommand{\bea}{\bed\begin{array}{rl}}
\newcommand{\eea}{\end{array}\eed}

\newcommand{\barray}{\begin{array}{ll}}
\newcommand{\earray}{\end{array}}
\newcommand{\lf}{\lfloor}

\newcommand{\Spvek}[2][r]{%
  \gdef\@VORNE{1}
  \left(\hskip-\arraycolsep%
    \begin{array}{#1}\vekSp@lten{#2}\end{array}%
  \hskip-\arraycolsep\right)}

\def\vekSp@lten#1{\xvekSp@lten#1;vekL@stLine;}
\def\vekL@stLine{vekL@stLine}
\def\xvekSp@lten#1;{\def\temp{#1}%
  \ifx\temp\vekL@stLine
  \else
    \ifnum\@VORNE=1\gdef\@VORNE{0}
    \else\@arraycr\fi%
    #1%
    \expandafter\xvekSp@lten
  \fi}

\def\sqr#1#2{{\vcenter{\vbox{\hrule height.#2pt
\hbox{\vrule width.#2pt height#1pt \kern#1pt
\vrule width.#2pt} \hrule height.#2pt}}}}

\newcommand{\Se}{\mathbb{S}}
\newcommand{\E}{\mathbb{E}}
\newcommand{\PP}{\mathbb{P}}
\newcommand{\RR}{\mathbb{R}}

\newcommand{\ba}{\bar a}

\newcommand{\bc}{\bar c}

\newcommand{\N}{\mathbb{N}}

\newcommand{\eproof}{\indent\vrule height6pt width4pt depth1pt\hfil\par\medbreak}

\def\a{\alpha} \def\ep{\epsilon} \def\uu{\upsilon}   \def\ka{\kappa}

\def\g{\gamma} \def\dis{\displaystyle}  
   
\def\va{\vartheta}    
  \def\r{\rho} \def\s{\sigma}

 \def\de{\delta}  \def\t0{\tau_0}
\def\ve{\vee}  \def\nn{\nonumber}

  \def\Gt{\mathcal{G}}
\def\T{\tau} \def\to{\tau_0} \def\G{\Gamma}  \def\D{\Delta}  
  
   \def\K{\times}

\def\1{\oslash} \def\2{\oplus} \def\3{\otimes} \def\4{\ominus}
\def\5{\circ} \def\6{\odot} \def\7{\backslash} \def\8{\infty}
\def\9{\bigcap} \def\0{\bigcup} \def\+{\pm} \def\-{\mp}
\def\la{\langle} 
\def\lf{\left} \def\rt{\right}\def\t{\triangle} 
\def\la{\label}\def\be{\begin{equation}}  \def\ee{\end{equation}}

\def\nn{\nonumber}

\def\bc{\begin{center}}       \def\ec{\end{center}}
\def\ba{\begin{array}}        \def\ea{\end{array}}
\def\be{\begin{equation}}     \def\ee{\end{equation}}
\def\bea{\begin{eqnarray}}    \def\eea{\end{eqnarray}}
\def\beaa{\begin{eqnarray*}}  \def\eeaa{\end{eqnarray*}}
\def\la{\label}

\begin{document}
\title{ Delay Feedback Control for  Switching Diffusion Systems
Based on  Discrete Time Observations\thanks{This work is entirely theoretical and the results can be reproduced using the methods described in this paper.}
}

\author{Xiaoyue Li\thanks{School of Mathematics and Statistics,
Northeast Normal University, Changchun, Jilin, 130024, China. Research
of this author  was supported by National Natural Science Foundation of China (11971096), the Natural Science Foundation of Jilin Province (No.20170101044JC), %the Education Department of Jilin Province (No.JJKH20170904KJ),
the Fundamental Research Funds for the Central Universities.}
\and Xuerong Mao \thanks{ Department of Mathematics and Statistics,
University of Strathclyde, Glasgow G1 1XH, U.K.
Corresponding author.  Research
of this author  was supported by
the Royal Society (WM160014, Royal Society Wolfson Research Merit Award),
the Royal Society and the Newton Fund (NA160317, Royal Society-Newton Advanced Fellowship) and the EPSRC (EP/K503174/1).}
\and Denis S. Mukama\thanks{School of Mathematics and Statistics,
Northeast Normal University, Changchun, Jilin, 130024, China.}
 \and Chenggui Yuan\thanks{Department of Mathematics, Swansea University, Bay Campus, SA1 8EN, UK.}}
\maketitle

\begin{abstract}
For the sake of saving time and costs the feedback control based on discrete-time observations is used to stabilize the switching diffusion systems. Response lags are required by most of physical systems and play a key role in the feedback control. The aim of  this paper is to design delay feedback control functions based on the discrete-time observations of the system states and the Markovian states in order for the controlled switching diffusion system (SDS) to be exponentially stable in $p$th moment and probability one as well as stable in $H_\8$. The designed control principles are implementable to stablize  quasi-linear and highly nonlinear SDSs. For quasi-linear SDSs the criteria are sharp that under the control with high strength the controlled SDSs will be stable (bounded) while under the weaker control they will be unstable (unbounded) in mean square. The sample and moment Lyapunov exponents are estimated which have close relationship with the time delays.
\end{abstract}

 {\it keywords} Brownian motion; Markov chain; Stochastic functional differential equations; Exponential stability; Moment boundedness; Lyapunov functional
%\end{keywords}
 {\it AMS}
\ 60H10;  93D15; 60J10
 % \end{AMS}
%\pagestyle{myheadings} \thispagestyle{plain} \markboth{LI, MAO, MUKAMA
%AND YUAN}{DELAY FEEDBACK CONTROL FOR SWITCHING DIFFUSION SYSTEMS}

\section{Introduction}\label{s1}

Switching diffusion systems (SDSs)  modulated by  Markov chains involving
   continuous dynamics and discrete events  provide more realistic models to describe the systems in many branches of science and industry which experience abrupt changes in their structures and parameters. Because of the wide range of applications,
dynamical properties of SDSs have been  investigated extensively
(see, e.g., \cite{ABG, my, yz09} and the references therein).   It is due to the Markovian switching that the dynamics of  SDSs may be drastically different from that of the systems without switching.  For example,  several counterexamples given in \cite{Pin1, Pin2} reveal  that the  recurrence or transience properties are opposite from their  subsystems' without switching.
For more properties such as the strong Feller, recurrence and stability please refer to \cite{shao0}, \cite{yf10} and the references therein.

One of the important issues in the study of SDSs is the automatic control, with consequent emphasis being placed on the their stabilization \cite{lu, mao13, mao08, mao14, qiu, shao, shao2, song, mao15}. Consider an unstable SDS  described by
\begin{equation} \label{q1}
dx(t) = f(x(t), r(t),t)dt + g(x(t), r(t),t)dB(t),
\end{equation}
where the state $x(t)$ takes values in $\RR^n$ and the mode $r(t)$ is a Markov chain taking values in a finite
space $\Se=\{1,2,\cdots, N\}$,  $B(t)$ is a Brownian motion.
 In order to stabilize this given system,  it is traditional  to
design a feedback control term $u(x(t),r(t),t)$
so that the controlled SDS (CSDS)
\begin{equation} \label{q2}
dx(t) = [f(x(t), r(t),t)+ u(x(t),r(t),t) ]dt + g(x(t), r(t),t)dB(t)
\end{equation}
 becomes stable.  Due to the requirement of the continuous-time observations for the state $x(t)$, it is difficult
 to implement  such a regular control.  In practice, very high frequent state observations are used instead of
 continuous-time observations and hence the control cost is expensive. For the sake of saving costs and easy operation Mao \cite{mao13} designed the feedback control based on the discrete-time observations (not necessarily high frequency), and developed the corresponding theory \cite{ds1,ds2,ds3} of  deterministic systems to stochastic versions. That is,  $u(x(\nu(t)),r(t),t) $ was designed, where $\nu(t):=[t/\T]\T$ with $\T >0$ being the duration between two consecutive observations, such that the controlled system
$$
 dx(t)= [f(x(t),r(t),t)+ u(x( \nu(t)),r(t),t) ]dt + g(x(t), r(t),t)dB(t),
$$
becomes stable in mean square. In the latter works \cite{mao14, mao15} much better lower bound on $\T$ was obtained while other types of asymptotic stability were studied.  However, from practical point of view it is sometimes necessary to design the  feedback control based on not only $x( \nu(t))$ but also  $r( \nu(t))$ (see, e.g., \cite{gg, song} for details).
Due to the continuity of $x(t)$ the deviation of $x(t)-x(\nu(t))$ may be small as long as $\T$ is sufficiently small. But the jump processes $r(\nu(t))$ and $r(t )$ may take different values in $\Se$ even if $\T$ is extremely small.  This problem was tackled by \cite{Li17,song}.  In particular, using different method from \cite{Li17,song},  Shao \cite{shao} obtained the stability in mean square for the linear controlled SDS
based on the discrete-time observations of both the system state $x(\cdot)$ and the Markov mode $r(\cdot)$. Shao and Xi \cite{shao2} went a further step to analyze the almost sure stability of the linear controlled SDS with the state-dependent  regime switching.

Response lags are often required by most physical systems, and play a crucial role in the feedback loops \cite{rm}.
Taking into account a time lag $\T_0$ ($>0$) between the time when the
observations for the state $(x( \nu(t) ), r( \nu(t) ))$ is made and the time when the feedback
control reaches the system, it is more
realistic to design the control dependent on the past discrete-time state pair $(x( \nu(t)-\to), r( \nu(t)-\to))$. To our best knowledge, the existing papers in the literature on stabilisation problems by delay feedback control are based on  the   observations of only system state $x(t-\to)$ or $x( \nu(t)-\to)$, for examples, \cite{lu, mao08, qiu, rm}.
Our main aim in this paper is to design the feedback control
 $u(x( \nu(t)-\to), r( \nu(t)-\to),t)$ ($\to\geq 0$) so that the delay controlled SDS (DCSDS)
\begin{equation} \label{q4}
 dx(t)= [f(x(t),r(t),t)+ u(x( \nu(t)-\T_0),r( \nu(t)-\T_0),t) ]dt + g(x(t), r(t),t)dB(t)
\end{equation}
 becomes stable in $p$th moment, with probability one or in $H_\8$.

 Mathematically speaking, this paper uses the strong ergodicity theory of Markov chains and the asymptotic analysis techniques of stochastic functional differential equations (SFDEs), which are completely different from those used in the
 papers \cite{mao13, mao14, qiu, shao, shao2, song, mao15} mentioned above.
 Various criteria on the uniform boundedness and different kinds of  stability  will be established for the DCSDS \eqref{q4} when their coefficients are either quasi-linear or highly nonlinear. The main contributions of this paper are highlighted as follows.
\begin{itemize}
\item For the quasi-linear DCSDSs \eqref{q4}, we give sharp criteria  on the uniform boundedness of the solution  in infinite horizon as well as  exponential stability in mean square.   That is, by a feedback control satisfying a proposed condition, the solution  will be uniformly bounded or exponentially stable in mean square, while it will be unbounded or unstable under a slightly weaker  control.
    The explicit rates of  the convergence and  divergence are obtained.

\item For the nonlinear DCSDSs \eqref{q4}, we give the criteria on the feedback control for the solution to be exponentially stable in $p$th moment and probability one as well as in $H_\8$. The sample and moment Lyapunov exponents are estimated, which describe the convergence speed that $x(t)$ tends to $0$ in $p$th moment and in sample path.

\item The lower bound on $\tau^*$ is obtained explicitly so that the feedback control will stablize the given system as long as $\T+\T_0 <\T^*$.
How the values of $\T$ and $\to$ affect  the Lyapunov exponents is also investigated.

\end{itemize}

The rest of the paper is organised as follows. Section \ref{s2} begins with notations and preliminaries on  the properties of the exact solutions. Section \ref{s3} focuses on the quasi-linear DCSDS \eqref{q4}. The sharp criteria on boundedness (unboundedness) and stability (unstability) are established. The convergence and divergence rates are estimated.
 Section \ref{s4} pays attention to the stability analysis for the highly nonlinear DCSDSs \eqref{q4}.
 Under the conditions on the existence of the global regular solution of \eqref{q1} and its boundedness in $p$th moment,
 it will be shown that the controlled system \eqref{q4}  preserves  the boundedness.
The  lower bounds on both $\T$ and $\to$ are also given explicitly.
 The control principles for the controlled system \eqref{q4} to be exponentially stable in $p$th moment or in probability one or in $H_\8$ are provided.   Furthermore, the sample and moment Lyapunov exponents are estimated.
in Section \ref{s5}, an example with computer simulations  is discussed to illustrate the theoretical results.

\section{Preliminary}\la{s2}
Throughout this paper,  we  use the following notations. If $A$ is a vector or matrix, its transpose is denoted by $A^T$ and its trace norm is denoted by $|A|=\sqrt{\mathrm{trace}(A^TA)}$.  For vectors  or matrixes $A$ and $B$ with compatible dimensions, $A B$ denotes the usual matrix multiplication.    For any sequence $ \{ c_i\}_{1\leq i\leq N}$ ($N\in \N$),
 define $ \hat{c}=\min_{1\leq i\leq N}   c_i$ and $ \check{c}=\max_{1\leq i\leq N} c_i $.  For any $a, b\in \mathbb{R}$,
 $a\vee b:=\max\{a,b\}$, and $a\wedge b:=\min\{a,b\}$.

  Let $( \Omega, \cal{F}, \PP )$ be a complete  probability space and $\mathbb{E}$ denote  the expectation with respect to $\PP$.
 Let $B(t)  = (B_1(t),\cdots,B_m(t))^T$
be an $m$-dimensional Brownian motion defined on the probability space.
Let $r(t)$, $t\ge 0$, be a right-continuous Markov chain on the probability space
taking values in a finite state space $\mathbb{S}=\{1, 2, \cdots, N\}$ ($N<\8$) with generator
$\G=(\g_{ij})_{N\K N}$ given by
\be\la{q-10}
\PP\{r(t+\D)=j | r(t)=i\}=
\begin{cases}
 \g_{ij}\D + o(\D) & \hbox{if }
i\not= j, \\
        1+\g_{ii}\D + o(\D) & \hbox{if }  i=j,
        \end{cases}
\ee
where $\D\downarrow 0$, $o(\D)$ means $\lim_{\D\rightarrow 0} o(\D)/\D=0$.  Here we assume $\G$ is {\it conservative }
 (i.e.
$
-\g_{ii}  =  \sum_{j\not= i} \g_{ij}, ~ \forall i\in\Se
$) and {\it irreducible} (i.e.
the linear equations $
\pi \Gamma=0$ and $  \sum_{i=1}^{N}\pi_i=1
$
has a unique solution $\pi=(\pi_1, \dots, \pi_N)\in {\mathbb R}^{1\times N}$  satisfying $\pi_i>0$ for each $i\in \mathbb{S}$). This solution is termed a stationary distribution. For a sequence $ \{ c_i\}_{1\leq i\leq N}$, we will often write $c(i)=c_i$ and set $c=(c_1,\cdots, c_N)^T$, define $\pi c = \sum_{i=1}^N \pi_i c_i$.  We assume that the Markov chain $r(\cdot)$ is independent of the Brownian motion
$B(\cdot)$.  Suppose ${\{ {\cal{F}}_{t}\}} _{t \geq 0}$ is a filtration defined on this probability space satisfying the usual conditions (i.e., it is right continuous and $\mathcal{F}_0$ contains all $\mathbb{P}$-null sets) such that $B(t)$ and $r(t)$
are $ {\cal{F}}_{t}  $ adapted. Denote by
 $\mathcal{G} $ the $\sigma$-algebra  generated by   $\{  r(t)\}_{0\leq t<\8}$.  We also denote the conditional  expectation $\E(\cdot|\mathcal{G})$ by  $\E_{\Gt}(\cdot)$.

In the paper, we use the feedback control function with a  simple form
$u(x,i,t)$ $=- \a(i) x $ for $(x,i,t)\in \RR^n\K \Se\K \RR_+$, where $\a(i)$'s are all nonnegative constants.
 Suppose that the underlying system is described by the DCSDS
(\ref{q4}) with the initial data
\begin{align}\label{q7}
   x(t)  =x_0\in \RR^n,  ~~~~   r(t)=i_0\in \Se,~  -\T_0\le t\le 0,
\end{align}
while the coefficient  functions
$
f:  \RR^n \K \Se\K \RR_+ \rightarrow \RR^n \quad\hbox{and}\quad
g:  \RR^n \K   \Se\K \RR_+ \rightarrow\RR^{n\K m}
$
 satisfy the {\it local Lipschitz} condition, namely, for any real number $R>0$, there exists a positive constant $K_R$ such that
$$
 |f(x,i,t)-f(\bar x,i,t)| \ve |g(x, i,t)-g(\bar x, i,t)| \le K_R |x-\bar x|
$$
 for all $x, \bar x \in \RR^n$ with $|x| \ve |\bar x| \le R$ and all $(i,t)\in \Se\K \RR_+$.
 It is well known that the local Lipschitz conditions of the coefficients  only guarantee that the SDS (\ref{q1}) has a unique maximal local solution, which  may explode to infinity at a finite time.   To avoid such a possible explosion, we  impose the following Khasminskii-type condition.

\begin{assp} \label{A1}
 Assume that there exist positive constants $ A ,  C$, and $p\geq 2 $ such that
$$
 x^T f(x, i,t) + \frac{p-1}{2} |g(x,i,t)|^2 \le C+A |x|^2,~~(x,i,t) \in \RR^n\K \Se \K \RR_+.
$$
  \end{assp}

 We prepare the regularity for the solutions of  SDS (\ref{q1}) and DCSDS (\ref{q4}), respectively, as follows.

\begin{lem} {\rm \cite[p. 93, Theorem 3.17]{my}}\label{l1}
Under Assumption \ref{A1}, the SDS (\ref{q1}) with the initial data $(x(0), r(0))=(x_0,i_0)\in \RR^n\times \Se$  has a unique global solution $x(t)$ on $[0,\8)$.
\end{lem}

In a similar way as \cite[p. 89, Theorem 3.13]{my} was proved, we can show:

\begin{lem} \label{l2}
Under Assumption \ref{A1}, the DCSDS (\ref{q4}) with the initial data \eqref{q7}  has a unique global solution $x(t)$ on $[0,\8)$.
\end{lem}

In \eqref{q4} the feedback control  depends on  the term $\a (r(\nu(t)-\to))$.  To analyze the asymptotic property we need a number of  new notations and recall some results from \cite{Ba}.
 For any vector $\mu=(\mu_1,\dots, \mu_N)^T$,  any constant $ l>0$,  define
\begin{align}\label{q-6}
\mathrm{diag}(\mu):= \mathrm{diag}(\mu_1,  \dots,   \mu_N),~\Gamma_{l,\mu}:= \Gamma-{l}\mathrm{diag}(\mu),~\eta_{l,\mu}:= -\!\!\!\max_{\lambda\in \mathrm{spec}(\Gamma_{l,\mu})}\mathrm{Re}(\lambda),
\end{align}
  where
$\mathrm{spec}(\Gamma_{l,\mu})$ and $\mathrm{Re}(\lambda)$ denote the spectrum of $\Gamma_{l,\mu}$ (i.e.  the multiset of its eigenvalues)  and the real part of $\lambda$, respectively.

  \begin{lem}\la{l5}{\rm \cite[Proposition 4.1, Proposition 4.2]{Ba}}
For any $l>0$, there are two positive constants
$K_1(l)$ and $K_2(l)$ such that  for any $t>0$
$$
K_1(l) e^{-\eta_{l,\mu}t }\leq \E \lf( e^{-l \int_0^{t} \mu(r(z)) dz}\rt) \leq K_2(l) e^{-\eta_{l,\mu}t }.$$
Moreover, if $
\pi \mu >0,
$   there is a constant $\kappa_\mu>0$ such that  $\eta_{l,\mu}>0$ for  $l\in (0,\kappa_\mu)$ but $\eta_{l,\mu}<0$ for $l>\kappa_\mu$.   Furthermore, if $\hat{\mu}  \geq  0$, $\kappa_\mu=\8$; if $\hat{\mu}  < 0$, $\kappa_\mu\in  (0,  \dis\min_{i\in \mathbb{S}, \mu_i<0}\left\{  \gamma_{ii}/\mu_i\right\})$.
  \end{lem}

In order to obtain the dynamical behaviors of the solutions of DCSDS \eqref{q4} we need to investigate the asymptotic properties of $\a (r(\nu(t)-\to))$. Firstly we redefine two Markov chains.
Let $n_0=[\to/\T]$, $\delta=(n_0+1)\T-\to$, $\tilde{r}(t):=r(t+(n_0+1)\T)$ for $t\geq 0$,  and  $ {r}_n:= r(n\T )$ for any integer $n\geq 0$.   Then $\{ {r}_n\}_{n\geq 0}$ is a skeleton process of Markov chain $\{r(t)\}_{t\geq 0}$, which is a discrete-time homogeneous Markov chain on $\Se$. Its transition probability matrix is $(P_{ij})_{N\times N}$ with $P_{ij}=\PP(r( \T+\de)=j|r(\de)=i )$. By virtue of Lemma \ref{l5} we can obtain the following results.

\begin{lem}\la{l4}
Let $h=(h_1,\cdots, h_N)^T$ such that $\pi \a<\pi h$.  Then, for  any  constants $l>0$ and  $0< \ep <\pi h-\pi \a$,  there is a constant $ T>0$ such that for any $s\in [0,\8)$
\begin{equation}\la{q9}
\E \lf( e^{l\int_s^{s+t}  (h(r(z))-\a(r(\nu(z)-\to))    )dz}\rt) \geq   K_1(l,\a-h)e^{   l (\pi h -\pi \a -\ep )t},~~~~t\geq T.
\end{equation}
 where $ K_1(l,\a-h):=e^{-l (2\T+\to) (\check{\a}+2\max_{i\in \Se} |h(i)| ) } $.
  \end{lem}

\noindent{\it Proof.}~~ To highlight the initial values, we let $\{ {r}^i(t)\}_{t\geq 0}$  and $\{ {r}^i_n\}_{n\geq 0}$ be the Markov chains starting from state $i\in\Se$ at $t=0$ and $n=0$, respectively. For any $i\in \Se$, since   $\{ {r}^i(t)\}_{t\geq 0}$ and $\{ {r}^i_n\}_{n\geq 0}$  are ergodic and has the same stationary distribution $(\pi_1, \dots, \pi_N)$, by the strong ergodic theorem and the boundedness of  $h(\cdot)$ and $ \a(\cdot )$,   we have
 \begin{equation}\la{q10}
 \begin{split}
 &\lim_{t\rightarrow\8} \frac{1}{t} \int_0^{ t}  (h( {r}^i(z +\to ))-\a( {r}^i(\nu(z) )) )dz\\
     &=   \pi h -\lim_{t\rightarrow\8} \frac{1}{t} \lf[   \T \sum_{j=0}^{[t/\T]-1 }    \a( {r}^i_j)+ \a( {r}^i_{[t/\T]} )(t-\nu(t))\rt]\\
 &=   \pi h -\lim_{t\rightarrow\8} \Big(\frac{\nu(t) }{t}\cdot   \frac{  \sum_{j=0}^{[t/\T] -1}    \a( {r}_j^i) }{[t/\T] } +\frac{ \a( {r}^i_{[t/\T]} )(t-\nu(t))}{t} \Big)\\
 &=   \pi h-\pi\a,~~~~\PP_2 -\hbox{a.s.}
 \end{split}
 \end{equation}
 By virtue of the Fatou lemma  (see, e.g. \cite[p.187, Theorem 2]{Shiryaev}),  for any given constant $l>0$, we have
 \begin{equation*}
 \liminf_{t\rightarrow \8}  \E \lf( \frac{l}{t}\int_0^{t}  (h( {r}^i(z+\to)) -\a( {r}^i(\nu(z) )) )dz  \rt) \geq  l(  \pi h-\pi\a).
\end{equation*}
If $\pi \a<\pi h$,   for any  $0<\ep<  \pi h-\pi \a, $ there is a constant $ T>0$ such that
$$
 \E \lf( \frac{l}{t}\int_0^{t}  (h( {r}^i(z+\to)) -\a( {r}^i(\nu(z) )) )dz  \rt) \geq  l(  \pi h-\pi\a-\ep)>0 ,~~~~i\in \Se,~ t\geq T
$$
holds. This implies
\begin{equation}\la{q11}
 \E \lf(  {l} \int_0^{t}  (h( {r}^i(z+\to)) -\a( {r}^i(\nu(z) )) )dz  \rt) \geq  l(  \pi h-\pi\a-\ep)t ,~~~~i\in \Se, ~t\geq T.
\end{equation}
Due to the Jensen inequality and the homogeneousness of Markov  chain   $\{ {r} (t)\}_{t\geq 0}$, we know that for any $s> 0$
\begin{equation}\la{q13}
\begin{split}
&  \E \lf( e^{ l\int_0^{t}  (h(\tilde{r} (z+s )) -\a(\tilde{r} (\nu(z+s )-\to)))dz } \rt)\\
& \geq    e^{  \E \lf(l\int_0^{t}  (h(\tilde{r} (z+s )) -\a(\tilde{r} (\nu(z+s )-\to)))dz  \rt)}\\
& = e^{ \E \lf( \E\lf(l\int_0^{t}   h(\tilde{r}  (z+s )) -\a(\tilde{r} (\nu(z+s )-\to))dz|\tilde{r} (\nu(s)-\to)   \rt) \rt) }\\
& = e^{  \sum_{j\in \Se}\E \lf(  I_{\{\tilde{r} (\nu(s)-\to)=j\}}\E\lf( l\int_0^{t}  h(\tilde{r} (z+s )) -\a(\tilde{r} (\nu(z+s )-\to))dz  |\tilde{r} (\nu(s)-\to)=j  \rt) \rt)}\\
 &=  e^{ \sum_{j\in \Se}\E \lf(  I_{\{ {r} (\nu(s) + \de)=j\}}\E\lf( l\int_0^{t}  h( {r}^j(z+ \de_s+\to )) -\a( {r}^j(\nu(z+ \de_s ) ))  dz  \rt) \rt) } ,
\end{split}
\end{equation}
where  $\de_s:=s-\nu(s)$  for any $s\geq 0$. From $0\leq \de_s<\T$,
one observes that for each $j\in \Se$
 \begin{equation}\la{q-1}
 \begin{split}
&\int_0^{t} (h( {r}^j(z+\de_s +\to) -\a( {r}^j(\nu(z+\de_s) )))dz \\& =\int_{\de_s}^{t+\de_s}  ( h( {r}^j(z+\to  ))-\a( {r}^j(\nu(z) )))dz\\
&=\int_{0}^{t}   ( h( {r}^j(z +\to ))-\a( {r}^j(\nu(z)  )))dz + \int_{t}^{t+\de_s}   ( h( {r}^j(z+\to  ))-\a( {r}^j(\nu(z) )))dz\\
&~~~~~~~~~~~- \int_{0}^{\de_s}  ( h( {r}^j(z +\to ))-\a( {r}^j(\nu(z) )))dz\\
& \geq \int_{0}^{t}   ( h( {r}^j(z+\to  ))-\a( {r}^j(\nu(z) )))dz- 2\T\max_{i\in \Se} |h(i)|-\T \check{\a}.
\end{split}
\end{equation}
Inserting \eqref{q-1} into \eqref{q13}, then using \eqref{q11}, we obtain that for any $i\in \Se$, $s\geq 0$
\begin{equation}\la{q-3}
\begin{split}
&  \E \lf( e^{ l\int_0^{t}  (h(\tilde{r} (z+s )) -\a(\tilde{r} (\nu(z+s )-\to)))dz } \rt)\\
&\geq e^{-\T l(\check{\a}+2\max_{i\in \Se} |h(i)| )}e^{\sum_{j\in \Se}\PP({r} (\nu(s) + \de)=j) \E  \lf(l\int_0^{t}  (h( {r}^j(z+\to  )) -\a( {r}^j(\nu(z ) )))dz  \rt) }\\
&\geq e^{-\T l(\check{\a}+2\max_{i\in \Se} |h(i)| ) } e^{    l (\pi h -\pi \a  -\ep )t},~~~~ t\geq T .
\end{split}
\end{equation}
This, together with the definition of $\tilde{r} (t)$, implies that for any $s\in [(n_0+1)\T,\8)$
\begin{equation}\la{q-4}
\E \lf( e^{l \int_s^{s+t}  (h(r  (z)) -\a(r (\nu(z)-\to))    )dz}\rt) \geq e^{-l\T (\check{\a}+2\max_{i\in \Se} |h(i)| ) } e^{   l (\pi h -\pi \a  -\ep )t},~~~~ t\geq T.
\end{equation}
By the similar way as \eqref{q-1} we know that for any $s\in [0, (n_0+1)\T)
$,
\be\la{q-2}\begin{split}
&\int_{s}^{s+t} ( h( {r} (z   ))-\a( {r} (\nu(z  ))-\to)))dz\\
&\geq \int_{(n_0+1)\T}^{t+(n_0+1)\T}   ( h( {r} (z  ))-\a( {r}  (\nu(z) -\to)))dz-(n_0+1)\T(\check{\a}+2\max_{i\in \Se} |h(i)|).
\end{split}\ee
This together with \eqref{q-4} implies that for any $s\in [0,\8)$
\begin{equation}\la{q-5}
\E \lf( e^{l \int_s^{s+t}  (h(r (z)) -\a(r (\nu(z)-\to))    )dz}\rt) \geq e^{-l\T (n_0+2) (\check{\a}+2\max_{i\in \Se} |h(i)| ) } e^{  l (\pi h-\pi \a  -\ep )t},~~  t\geq T.
\end{equation} The required assertion \eqref{q9} follows.
\eproof

\begin{lem}\la{l6} Let $h=(h_1,\cdots, h_N)^T$ such that  $\pi \a>\pi h $. Then, for  any  constant  $ 0<l<\ka_{\a-h}$, if $\T< \bar{\T}(l,\a-h)$, there are positive  constants $K_3 ({l,\a-h})$ and $\zeta_{l,\a-h}^\T $ defined by \eqref{q-18+} such that for any $t\geq 0$
\begin{equation}\la{q-8}
\E \lf( e^{l\int_0^{t}  (h({r}(z+\to))-\a( {r}(\nu(z)))    )dz}\rt) \leq K_3 ({l,\a-h}) e^{- \zeta_{l,\a-h}^\T  t},
\end{equation} where   $\bar\T=\bar{\T}(l,\a-h)$ is the solution of the equation (in $\T$)
 \be\la{q-18}\ep \max_{j\in \Se}\{ - \g_{jj}\}   (e^{   \frac{\T l \check{\a }(1+\ep) }{\ep} }-1 )=\eta_{l(1+\ep), \a-h},~~~~\ep:=[(\ka_{\a-h}-l)/2l]\wedge 1.\ee
  \end{lem}

\noindent{\it Proof.}~~ One observes
\begin{align}\la{+q-15}
  \int_{0}^{t} h( {r}(z+\to))dz
 =\int_{ \to}^{t+\to} h(r(z)) dz  \leq \int_{0}^{t}h(r(z)) dz+
 2\to \max_{i\in \Se} |h(i)| .\end{align} By H$\ddot{\hbox{o}}$lder's inequality, we obtain that for any $ 0<l<\ka_{\a-h}$ and the given $\ep>0$,
 \begin{equation}\la{q-9}
 \begin{split}
&\E \lf( e^{l\int_0^{t}  (h( {r}(z+\to))-\a( {r}(\nu(z)))    )dz}\rt) \\
&\leq e^{ 2l\to  \max_{i\in \Se} |h(i)| }\E \lf( e^{l\int_0^{t}  (h( {r}(z ))-\a( {r}(\nu(z)))    )dz}\rt) \\
&\leq e^{ 2l\to  \max_{i\in \Se} |h(i)| } \Big( \E e^{  (1+\ep)  l\int_0^{t}  (h( {r}(z ))-\a( {r}(z ))    )dz}\Big)^{\frac{1}{1+\ep}}\\\
&~~~~~~~~~~~~~~~~~~~~~~~~~~~~~\times
\Big( \E  e^{   \frac{l(1+\ep)}{\ep}   \int_0^{t}  (\a( {r}(z ))-\a( {r}(\nu(z)))    )dz}\Big)^{\frac{\ep}{1+\ep}}.
 \end{split}
\end{equation}
By virtue of Lemma \ref{l5} we know that
\begin{align}\la{q-13}
& \Big( \E e^{  (1+\ep)  l\int_0^{t}  (h( {r}(z ))-\a( {r}(z ))    )dz}\Big)^{\frac{1}{1+\ep}} \leq [K_2((1+\ep)l) ]^{\frac{1}{1+\ep}}e^{-\frac{\eta_{l(1+\ep), \a-h}}{1+\ep}t}.
\end{align}
 On the other hand
\begin{align}\la{q-11}
    \E  e^{   \frac{l(1+\ep)}{\ep}   \int_0^{t}  ( \a( {r}(z ))-\a( {r}(\nu(z) )))    dz}
     &\leq      \E  e^{   \frac{l(1+\ep)}{\ep}  \sum_{i=0}^{[t/\T]} \int_{i\T}^{(i+1)\T}  |\a( {r}(z ))-\a( {r}(i\T ))   |dz } \nn\\
     &=    \E \Big(\prod_{i=0}^{[t/\T]}  e^{   \frac{l(1+\ep)}{\ep} \int_{i\T}^{(i+1)\T}  |\a( {r}(z ))-\a( {r}(i\T ))   |dz }\Big).
 \end{align}
For any nonnegative integer $i$, the Jensen inequality shows
\begin{align*}
 &\E(  e^{   \frac{l(1+\ep)}{\ep} \int_{i\T}^{(i+1)\T}   |\a( {r}(z))-\a( {r}(\nu(z) ))   |dz}| {r}(i\T))\\
&\leq  \E( \frac{1}{\T} \int_{i\T}^{(i+1)\T} e^{   \frac{\T l(1+\ep) }{\ep} |\a( {r}(z ))-\a( {r}(\nu(z)))   | }dz| {r}(i\T)).
 \end{align*}
It is known that the waiting time for the next jump of the Markov chain $ {r}(\cdot)$ from current state $j$ obeys the exponential distribution with parameter $-\g_{jj}$ (see, e.g., \cite[p. 16, Proposition 2.8]{An}). Thus, noticing that $1-e^{-x}\leq x$ for $x\geq 0$, we have
 \begin{align*}
 &\E(  e^{   \frac{l(1+\ep)}{\ep} \int_{i\T}^{(i+1)\T} |\a( {r}(z))-\a( {r}(i\T ))   |dz}| {r}(i\T) )\nn\\
  &\leq \frac{1}{\T} \sum_{j=1}^{N}I_{\{ {r}(i\T)=j\}}\int_{i\T}^{(i+1)\T} \E  \Big( I_{\{ {r}(v)=j,~\forall v\in [i\T ,(i+1)\T] \}}\nn \\
    &~~~~~~+    I_{\{\exists v\in [i\T ,(i+1)\T],~  {r}(v)\neq j\}} e^{   \frac{\T l \check{\a } (1+\ep) }{\ep}}| {r}(i\T)=j\Big)dz\nn\\
    &\leq  \frac{1}{\T} \sum_{j=1}^{N}I_{\{ {r}(i\T)=j\}}\int_{i\T}^{(i+1)\T}\Big(  \PP  (  {r}(v)=j,~\forall v\in[i\T ,(i+1)\T]  | {r}(i\T)=j)\nn\\
    &~~~~~+  e^{   \frac{\T l \check{\a }(1+\ep) }{\ep} }\PP  (\exists v\in [i\T ,(i+1)\T],~  {r}(v)\neq j | {r}(i\T)=j)  \Big)dz\nn\\
     &=  \frac{1}{\T} \sum_{j=1}^{N}I_{\{ {r}(i\T)=j\}}\int_{i\T}^{(i+1)\T}\Big(  e^{\g_{jj}\T} + e^{   \frac{\T l \check{\a }(1+\ep) }{\ep} }(1- e^{\g_{jj}\T}) \Big)dz\nn\\
     &=    \sum_{j=1}^{N}I_{\{ {r}(i\T)=j\}} \Big(  1 + (e^{   \frac{\T l \check{\a }(1+\ep) }{\ep} }-1)(1- e^{\g_{jj}\T}) \Big)\leq   1+  \T \Lambda_\T(l),
 \end{align*}
where \be\la{q-16}\Lambda_\T(l):=\max_{j\in \Se}\{ - \g_{jj}\}   (e^{   \frac{\T l \check{\a }(1+\ep) }{\ep} }-1 ).
\ee Inserting this inequality into \eqref{q-11} yields that
\begin{equation}\la{q-12}
\begin{split}
   &\E  e^{   \frac{l(1+\ep)}{\ep}   \int_0^{t}   \a( {r}(z ))-\a( {r}(\nu(z) ))    dz}\\
    &\leq    \E \Big[\E\Big(\prod_{i=0}^{[t/\T]
    }  e^{   \frac{l(1+\ep)}{\ep} \int_{i\T}^{(i+1)\T}  |\a( {r}(z))-\a( {r}(i\T ))   |dz }| {r}([t/\T]\T)\Big)\Big]\\
   &\leq       \E \Big[\prod_{i=0}^{[t/\T]-1
    }  e^{   \frac{l(1+\ep)}{\ep} \int_{i\T}^{(i+1)\T}  |\a( {r}(z))-\a( {r}(i\T ))   |dz }\\
 &~~~~~~~~~~~~~~~~~~~~~~~~~~~~   \times \E\Big(   e^{   \frac{l(1+\ep)}{\ep} \int_{[t/\T]\T}^{([t/\T]+1)\T}  |\a( {r}(z))-\a( {r}([t/\T]\T+\de))   |dz }| {r}([t/\T]\T)\Big)\Big]\\
 &\leq
 \Big(1+      \T \Lambda_\T(l) \Big)^{[t/\T]+1}
 \leq
 e^{([t/\T]+1)   \T  \Lambda_\T(l)}
   \leq
 e^{(t +\T) \Lambda_\T(l)}.
  \end{split}
\end{equation}
  Thus, inserting \eqref{q-13} and \eqref{q-12} into \eqref{q-9} we obtain
   \begin{equation}\la{q-17}
 \begin{split}
 \E \lf( e^{l\int_0^{t}  (h( {r}(z+\to))-\a( {r}(\nu(z) ))    )dz}\rt)
 \leq K_3({l,\a-h}) e^{- \zeta_{l,\a-h}^\T t},
 \end{split}
\end{equation} where \be\la{q-18+}\begin{split}
K_3({l,\a-h}):&=[K_2((1+\ep)l)]^{\frac{1}{1+\ep}} e^{\frac{  \ep \T \Lambda_\T(l)}{1+\ep}+ 2l\to  \max_{i\in \Se} |h(i)| }, \\
\zeta_{l,\a-h}^\T:&=\frac{\eta_{l(1+\ep), \a-h}- \ep  \Lambda_\T(l)}{1+\ep}.\end{split}\ee
Using $\T<\bar{\T}(l,\a-h) $ and the definition of $\Lambda_\T(l)$, we know that $\zeta_{l,\a-h}^{\T}>0$. Therefore,  the required assertion follows.\eproof

\begin{lem}\la{l7} Let $h=(h_1,\cdots, h_N)^T$ such that $\pi \a>\pi h $. For  any  constant  $ 0<l<\ka_{\a-h}$, if $\T< \bar{\T}(l,\a-h)$, there is a positive  constant  $K_2  ({l,\a-h})$  such that for any $s\geq 0$
\begin{equation}\la{q-20}
\E \lf( e^{l\int_s^{s+t}  (h({r}(z))-\a({r}(\nu(z)-\to))    )dz}\rt) \leq K_2 ({l,\a-h}) e^{- \zeta_{l,\a-h}^\T  t}, ~~t\geq 0,
\end{equation} where  $K_2  ({l,\a-h}):=K_3  ({l,\a-h})e^{ l (2\T+\to) (\check{\a}+2\max_{i\in \Se} |h(i)| ) }$,  $\bar{\T}(l,\a-h)$  and $\zeta_{l,\a-h}^\T $ are given in Lemma \ref{l6}.
  \end{lem}

\noindent{\it Proof.}~~ Using the homogeneousness of $ {{r}}(\cdot)$ and the definition of $ {{\tilde{r}}}(\cdot)$, we know that for any $t\geq 0$, $s> 0$,
\begin{align}\la{q+13}
&  \E \lf( e^{ l\int_0^{t}  (h(\tilde{r} (z+s )) -\a(\tilde{r} (\nu(z+s )-\to)))dz } \rt)\nn\\
& =   \sum_{j\in \Se}\E \lf(  I_{\{\tilde{r} (\nu(s)-\to)=j\}}\E\lf( e^{ l\int_0^{t}  (h(\tilde{r} (z+s )) -\a(\tilde{r} (\nu(z+s )-\to)))dz }  |\tilde{r} (\nu(s)-\to)=j  \rt) \rt) \nn \\
& =   \sum_{j\in \Se}\E \lf(  I_{\{{r} (\nu(s)+\de)=j\}}\E\lf( e^{ l\int_0^{t} (  h({r}^j(z+ \de_s+\to )) -\a({r}^j(\nu(z+ \de_s ))) ) dz } \rt) \rt).
\end{align}
 From $0\leq \de_s<\T$, by the similar way as \eqref{q-1},
one observes that for each $j\in \Se$
 \begin{equation}\la{+q-1}
 \begin{split}
&\int_0^{t} (h({r}^j(z+\de_s+\to ) -\a({r}^j(\nu(z+\de_s))))dz \\
 & \leq \int_{0}^{t}   ( h({r}^j(z+\to  ))-\a({r}^j(\nu(z) )))dz+ 2\T\max_{i\in \Se} |h(i)|+\T \check{\a}.
 \end{split}
\end{equation}
Inserting \eqref{+q-1} into \eqref{q+13}, then using  the fact that the estimate of \eqref{q-20} is independent of the initial value $ {r}(0 )$, we obtain that for any $i\in \Se$, $s\geq 0$,
\begin{align}\la{+q-3}
&  \E \lf( e^{ l\int_0^{t}  (h(\tilde{r} (z+s )) -\a(\tilde{r} (\nu(z+s )-\to)))dz } \rt)
 \leq  K_3  ({l,\a-h})e^{ l \T( \check{\a}+2  \max_{i\in \Se} |h(i)| ) }e^{- \zeta_{l,\a-h}^\T  t}.
\end{align}
This, together with the definition of $\tilde{r} (t)$, implies that for any $s\in [(n_0+1)\T,\8)$
\begin{equation}\la{+q-4}
\E \lf( e^{l \int_s^{s+t}  (h(r  (z)) -\a(r (\nu(z)-\to))    )dz}\rt)  \leq K_3  ({l,\a-h})e^{ l \T(  \check{\a}+2  \max_{i\in \Se} |h(i)| ) }e^{- \zeta_{l,\a-h}^\T  t}.
\end{equation}
By the similar way as \eqref{q-2} we know that for any $s\in [0, (n_0+1)\T)$,
\be\la{+q-2}\begin{split}
&\int_{s}^{s+t} ( h( {r} (z   ))-\a( {r} (\nu(z  ))-\to)))dz\\
&\leq \int_{(n_0+1)\T}^{t+(n_0+1)\T}   ( h( {r} (z  ))-\a( {r}  (\nu(z) -\to)))dz+(n_0+1)\T(\check{\a}+2\max_{i\in \Se} |h(i)|).
\end{split}\ee
This together with \eqref{+q-4} implies that for any $s\in [0,\8)$,
\begin{equation}\la{+q-5}
\E \lf( e^{l \int_s^{s+t}  (h(r (z)) -\a(r (\nu(z)-\to))    )dz}\rt)  \leq K_2  ({l,\a-h}) e^{- \zeta_{l,\a-h}^\T  t}.
\end{equation} The required assertion \eqref{q-20} follows.
\eproof

In order to deal with the asymptotic properties of the SFDE \eqref{q4}, we begin with defining two segments $\bar{x}_t(s):= \{x(t+s): -2(\T+\to) \leq s \leq 0\} $ and  $\bar{r}_t(s):= \{r(t+s): -2(\T+\to)\leq s \leq 0\} $ for $t\geq 0$. In order for $\bar{x}_t$ and $\bar{r}_t$ to be well defined on $0\leq t<2(\T+\to)$, we let $x(s)=x_0$ and $r(s)=i_0$ for $s\in [-2 (\T+\to), -\T_0)$. Moreover,  we enlarge the definition domains of $f$, $g$  and $u$. For any $(x,i,t)\in \RR^n\times \Se\times [-2(\T+\to), 0)$, let $f(x,i,t)=f(x,i,0)$, $g(x,i,t)=g(x,i,0)$, $u(x,i,t)=u(x,i,0)$.
In order to control the derivation from time delay in mean square, i.e. the value of $\E_{\Gt} |x(t)- x(\nu(t)-\to)|^2$,
 we define an auxiliary functional
  \begin{align}\la{q14}
 I(\bar{x}_t,\bar{r}_t, t)
 &:=\int_{-(\T+\to)}^0 \int_{t+s}^t  \lf[(\T+\to)|f(x(z),r(z), z) -\a(r(\nu(z)-\to))x(\nu(z)-\to)|^2 \rt.\nn \\
 &~~~~~~~~~~~~~~~~~~~~\lf.+ |g(x(z),r(z), z) |^2\rt]dzds.
 \end{align} For simplicity  we let $I( t)=I(\bar{x}_t,\bar{r}_t, t)$.
  A direct calculation arrives at \be\la{q15} dI(t)=J_1(t)dt-J_2(t)dt, \ee where
 \begin{align}
 J_1(t)&= (\T+\to)   \lf[(\T+\to)|f(x(t),r(t), t) -\a(r(\nu(t)-\to))x(\nu(t)-\to)|^2\rt.\nn \\
  &~~~~~~~~~~~~~~~~~~~~~~~~~ \lf.+ |g(x(t),r(t), t) |^2\rt],\la{q16}\\
  J_2(t)&= \int_{-(\T+\to)}^0    \lf[(\T+\to) |f(x(t+s),r(t+s), t+s)  \rt.\nn \\
 &  \lf.-\a(r(\nu(t+s)-\to))x(\nu(t+s)-\to)|^2+ |g(x(t+s),r(t+s), t+s) |^2\rt]ds. \la{q17}
\end{align}
By changing the integration order, we  get
 \be\la{q18} I(t)\leq  (\T+\to)J_2(t).\ee
Using the H$\ddot{\hbox{o}}$lder inequality and the It$\hat{\hbox{o}}$ isometry formula we go a further step to obtain
\be\la{q19}\begin{split}
&\E_{\Gt} |x(t)- x(\nu(t)-\to)|^2 \\
&= \E_{\Gt}|\int_{\nu(t)-\to}^t [f(x( s),r( s),  s) -\a(r(\nu(s)-\to))x(\nu(s)-\to)]ds\\
 & ~~~~~~~~~~~~ +g(x( s),r( s),  s)dB(s)  |^2 \\
 &\leq 2 \E_{\Gt} \int_{t-(\T+\to) }^t  (\T+\to) |f(x( s),r( s),  s) -\a(r(\nu(s)-\to))x(\nu(s)-\to)|^2 ds \\
 & ~~~~~~~~~~~~ +2\E_{\Gt} \int_{t-(\T+\to) }^t |g(x( s),r( s),  s)|^2ds =2\E_{\Gt} [J_2(t)].
\end{split}\ee

\section{Control of Quasi-linear Systems}\la{s3}
This section pays attention to design the control functions for the  solutions of quasi-linear SDSs    to be bounded in the infinite time horizon, and exponentially stable in $p$th moment and in $\PP-$a.s.

\subsection{Boundedness Control}
As is well-known that the unique solution of a linear SDS exists globally on $[0, \8)$, and its $p$th moment is finite on any finite time interval $[0, T]$. However its $p$th moment may be unbounded in infinite horizon $[0,\8)$. So it is necessary to design the control function $u(x(\nu(t)-\to), r(\nu(t)-\to), t)$ such that the solution of the  controlled system \eqref{q4} is bounded  in mean square in  $[0,\8)$.  To be precise we state the hypothesis of the linear growth condition on the coefficients.

 \begin{assp} \label{A2}
There exist positive constants  $\bar{K} $,  $ D_i $ and $ E_i  $  such that
     \begin{equation}\label{q20}
 |f(x, i,t)|\vee    |g(x, i,t)| \le \bar{K}  (1+ |x|) ,
 \end{equation} and
 \begin{equation}\label{q21}
   x^T f(x, i,t) + \frac{  1}{2} |g(x,i,t)|^2 \le E_i + D_i |x|^2
 \end{equation}
 hold for all $(x,i,t) \in \RR^n\K \Se \K \RR_+$.
 \end{assp}

\begin{thm}\la{t1}
Let   Assumption \ref{A2} hold and set $D=(D_1,\cdots, D_N)^T$.   Assume that $ \pi\a>\pi  D$ and  $ \ka_{\a-D}>2 $. If   $\T\leq {\tilde{\T}} := \bar{\T}(2,\a-D)/2$, $\T+\to\leq \T^*_1:=y_1\wedge y_2$, then the solution of DCSDS \eqref{q4} with  the initial condition \eqref{q7} has the property that \be\la{q23}
\sup_{0\leq t< \8} \E |x(t)|^2  <\8, \ee
where $\bar{\T}(\cdot, \cdot) $ is given as the solution of equation \eqref{q-18}, $y_1$ and $y_2$ are the  positive solutions of equations
$$\bar{\beta}_1(y)=\frac{\zeta^2}{2(8\check{\a}^2+\zeta^2)},~~~~\bar{\beta}_2(y)=\frac{\zeta^2}{ 8\check{\a}^2+\zeta^2 },$$ respectively. Here we write $\zeta =\zeta_{2,\a-h}^{{\tilde{\T}}}$ for short, $\bar{\beta}_1(\cdot)$, $\bar{\beta}_2(\cdot)$ are defined by \eqref{q29} and \eqref{q36} below, respectively.
\end{thm}
\noindent{\it Proof.}
Fix $0<\T\leq {\tilde{\T}} $ and $0<\T+\to\leq \T^*_1$.
 Using  \eqref{q20}  we compute $I(t) $ and $J_1(t)$ defined by \eqref{q14} and \eqref{q16}
\begin{align} \la{q24}
 J_1(t)
  &\leq2(\T+\to)  \lf[  \lf(2(\T+\to) ( \bar{K}^2+  \check{\a}^2)+  \bar{K}^2\rt) |x(t)|^2 \rt.\nn\\
  &~~~~~~~~\lf.+2\check{\a}^2(\T+\to) |x(t)-x(\nu(t)-\to) |^2 +   2\bar{K}^2 (\T+\to)  +  \bar{K}^2  \rt].
\end{align}
Inserting  \eqref{q24} into \eqref{q15} yields
\begin{align} \la{q25}
d I(t)
  \leq  \Big( {\bar{\beta}}_1(\T+\to)    |x(t)|^2    +\iota(\T+\to)    |x(t)- x(\nu(t)-\to)|^2 +{\bar{\beta}}_1(\T+\to)-J_2(t)\Big)dt,
\end{align}  where for any $y\geq 0$,
 \begin{align}    \bar{\beta}_1(y)   : = 2 y \lf[ 2y(  \bar{K}^2+ \check{\a}^2)+ \bar{K}^2\rt], ~~~~ \iota(y)  :=4 \check{\a}^2y^2.\la{q29}
   \end{align}
 Using the It$\hat{\hbox{o}}$ formula and the elementary inequality, by \eqref{q21}, for $\zeta>0$, we derive
  \begin{align}
 d|x(t) |^{2}
&\leq  \Big[ 2\check{ E}+  \Big(2D(r(t))  - 2 \a(r(\nu(t)-\to))+\frac{\zeta}{2}\Big)|x(t)|^2 \nn\\
&~~~~~~~~~~~~+ \frac{ 2\check{\a}^2}{\zeta} |x(t)- x(\nu(t)-\to)|^2\Big]dt+  2 x^T(t)  g(x(t),r(t),t)  dB(t).\la{q26}
 \end{align}
Define
$V(\bar{x}_t,\bar{r}_t, t)=(|x(t)|^2+\eta I(t))e^{-2  \int_0^t \phi(s) ds },$ % \ee
 where
 $\phi(s):= D(r(s)) -\a(r(\nu(s)-\to))+ {3\zeta }/{8},~ \eta:= {\zeta}/{2}+ {4\check{\a}^2}/{\zeta} .$
Since $\phi(s)$ has only a finite number of jumps in any  finite interval $[0, t]$, $\int_0^t  \phi(s)ds$ is differentiable. It follows from \eqref{q25} and \eqref{q26} that for any $t\geq 0$,
\begin{align} \la{q31}
d V(\bar{x}_t,\bar{r}_t, t)&\leq e^{-2  \int_0^t \phi(s) ds } \Big[-2 \eta \phi(t) I(t)  -(\frac{\zeta}{4}-\eta\bar{\beta}_1(\T+\to)  ) |x(s)|^2   \nn\\
&  ~~    + \lf(  \frac{2\check{\a }^2}{\zeta}+\eta\iota (\T+\to) \rt)  |x(s)- x(\nu(s)-\to)|^2 +(2\check{E}+ \eta \bar{\beta} _1(\T+\to)) \nn\\
  &~ ~ ~~ -  {\eta} J_2(t)\Big]dt+ 2 e^{-2  \int_0^t \phi(s) ds } x^T(t)  g(x(t),r(t),t)  dB(t) .
\end{align}
 Due to the increasing property of  $\bar{\beta}_1(y)$ in $y>0$,  we see
$\eta \bar{\beta}_1(\T+\to) \leq \zeta/4.$
 One observes from \eqref{q18} that
   \be\la{q33} -2 \eta \phi(t) I(t) \leq 2 \eta \check{\a} I(t)\leq 2 \eta\check{\a}  (\T+\to)J_2(t).\ee
    These, together with \eqref{q31}, imply
\begin{equation} \la{q34}
\begin{split}
d V(\bar{x}_t,\bar{r}_t, t)&\leq e^{-2  \int_0^t \phi(s) ds } \Big[
        \lf(  \frac{2\check{\a }^2}{\zeta}+\eta\iota (\T+\to) \rt)  |x(s)- x(\nu(s)-\to)|^2 \\
        &~ ~ ~~    +\lf(2\check{E}+ \frac{\zeta}{4}\rt)- ({\eta}-2{\eta}\check{\a}  (\T+\to)  )J_2(t)\Big]dt \\
  &~ ~ ~~ + 2 e^{-2  \int_0^t \phi(s) ds } x^T(t)  g(x(t),r(t),t)  dB(t).
  \end{split}
\end{equation}
Integrating \eqref{q34} on both sides and then taking the conditional expectation with respect to $\Gt$ and using \eqref{q19}, we  arrive  at
\begin{equation}
\begin{split}
 & \E_{\Gt}\lf(e^{-2  \int_0^t \phi(s) ds }|x(t)|^2\rt)+  \eta\E_{\Gt}  \lf(e^{-2  \int_0^t \phi(s) ds }I(t) \rt)\\
&\leq  |x_0|^2+\eta  I(0)  +\lf(2\check{E}+ \frac{\zeta}{4}\rt)\int_0^t e^{-2  \int_0^s\phi(z) dz  }ds \\&~~~~~~~~~~~~~~~~-  \lf(  {\eta}  -\frac{4\check{\a }^2}{\zeta}- \eta {\bar{\beta}}_2(\T+\to) \rt)   \int_0^t e^{-2  \int_0^s \phi(z) dz } \E_{\Gt} [J_2(s) ] ds,\la{q35}
\end{split}
\end{equation}
where for any constant $y\geq 0$, \be\la{q36}  {\bar{\beta}}_2(y)  :=2 \iota(y) +2  y \check{\a} = 2\check{\a}y( 4 \check{\a} y  +1 ). \ee
Due to the increasing property of  $\bar{\beta}_2(y)$ in $y>0$ as well as by the definition of $\eta$,
we see that
$    \eta \bar{\beta}_2(\T+\to)\leq  \zeta/2=\eta-4\check{\a }^2/\zeta.$
This  together with \eqref{q35} implies that
\begin{align*}  \E_{\Gt}|x(t)|^2 &\leq (|x_0|^2+ \eta  I(0)) e^{2  \int_0^t \phi(s) ds }  +\lf(2\check{E}+ \frac{\zeta}{4}\rt)\int_0^t e^{  2 \int_s^t\phi(z) dz  }ds.
\end{align*}
Taking expectation
on both sides yields
\begin{align}\la{q37} \E |x(t)|^2  &\leq  \lf(|x_0|^2+  \eta  I(0) \rt) \E\lf(e^{2  \int_0^t  \phi(s) ds}\rt)     +\lf(2\check{E}+ \frac{\zeta}{4}\rt)\int_0^t \E\lf(e^{2  \int_s^t  \phi(z)dz}\rt)ds   .
\end{align}
Since $ \ka_{\a-D}>2 $, by virtue of Lemma \ref{l7},  we have
 \begin{equation*}
\E \lf( e^{2 \int_s^t  (D(r(z)))-\a(r(\nu(z)-\to))    )dz}\rt) \leq K_2(2, \a-D)e^{ -    { \zeta }(t-s) },~~~~ t\geq  s,~s\geq 0,
\end{equation*}
which implies
\be\la{q38}\E\lf(e^{2  \int_s^t  \phi(z)dz}\rt) \leq K_2(2, \a-D)e^{ -    {\zeta }(t-s)/4 },~~~~ t\geq  s,~s\geq 0. \ee
It follows from \eqref{q38} that for any $t\geq 0$,
\begin{align}\la{q40}
  \int_0^t \E\lf(e^{2  \int_s^t  \phi(z)dz}\rt)ds
 \leq  K_2(2, \a-D) \int_{0}^{t } e^{ -    \frac{\zeta }{4}(t-s)  }ds  \leq  \frac{4}{\zeta}K_2(2, \a-D) .
 \end{align}
Inserting \eqref{q38} and \eqref{q40} into \eqref{q37}  arrives at
\begin{align*} \E |x(t)|^2
 &\leq   \lf(|x_0|^2+  \eta  I(0) \rt) K_2(2, \a-D) e^{ -    \frac{\zeta }{4}t  } + \lf(\frac{8\check{E}}{\zeta}+ 1\rt) K_2(2, \a-D)
\end{align*}
for $ t \geq 0 $. Thus, the required assertion \eqref{q23} follows. \eproof

Next we consider  the opposite aspect, namely, if the control strength is taken smaller value what will happen.
We investigate the longtime behavior of the mean square of the DCSDS \eqref{q4} in this case.

\begin{assp} \label{A3}
   Assume that there exist positive constants  $\bar{K }>0$,    and  $ d_i,~e_i $     such that
      \eqref{q20}  and
 \begin{equation}\label{q41}
 x^T f(x, i,t) + \frac{  1}{2} |g(x,i,t)|^2 \geq d_i |x|^2+e_i
 \end{equation}
 hold for all $(x,i,t) \in \RR^n\K \Se \K \RR_+$.
 \end{assp}

\begin{thm}\la{t2}
Let   Assumption \ref{A3} hold and assume that $\upsilon:=\pi  d -\pi\a >0$, where $d=(d_1,\cdots, d_N)^T$.
If $ 0<\T+\to < \T_*^1=y_3\wedge y_4\wedge y_5$,
 then   the solution of DCSDS \eqref{q4} with  the initial solution \eqref{q7} has the property that
 \be\la{q43}
\lim_{t\rightarrow \8} \E |x(t)|^2=\8,
\ee where $y_i$ $(i=3,4,5)$ are the maximum positive solutions of $$\bar{\beta}_1(y)=\frac{\uu(\hat{e}\wedge (\uu/2))}{ 2\check{\a}^2+\uu^2 },~~ \bar{\beta}_3(y)=\frac{\uu^2}{ 2\check{\a}^2+\uu^2 },~~\bar{\beta}_4(y)=\frac{\uu( {\hat{e}}/(2\check{d})+ |x_0|^2)}{ 2\check{\a}^2+\uu^2 },$$ respectively, $\bar{\beta}_1(\cdot)$, $\bar{\beta}_3(\cdot)$, $\bar{\beta}_4(\cdot)$ are defined by \eqref{q29}, \eqref{q57}, \eqref{q58}.
\end{thm}
\noindent{\it Proof.}
Fix  $ 0<\T+\to  \leq \T_*^1 $. Using  the elementary inequality and  \eqref{q41}, we derive    \begin{align}
 d|x(t) |^{2}
&\geq \Big[  (2d(r(t))  - 2 \a(r(\nu(t)-\to))-\uu)|x(t)|^2 - \frac{\a^2}{\uu} |x(t)- x(\nu(t)-\to)|^2  \nn\\
&~~~~~~~~~~~~+2\hat{e}\Big]dt+  2 x^T(t)  g(x(t),r(t),t)  dB(t).\la{q44}
 \end{align}
 Define
 ${U}(\bar{x}_t,\bar{r}_t, t)=\lf( { \hat{e}}/({2\check{d}})+|x(t)|^2-\eta I(t)\rt)e^{-2  \int_0^t \psi(s) ds },$  where
  $\psi(s):= d(r(s))$ $ -\a(r(\nu(s)-\to))-3\uu/4,~ \eta=\uu+ {2\check{\a}^2}/{\uu},~ \forall s\geq 0.$
 One notices that
 $ - ({ \hat{e}}/{ \check{d}} )\psi(t)\geq -\hat{e}$. This together with
 \eqref{q25} and \eqref{q44} implies that for any $t\geq 0$
\begin{equation} \la{q53}
\begin{split}
d U(\bar{x}_t,\bar{r}_t, t)&\geq e^{-2  \int_0^t \psi(s) ds } \Big[ \Big( 2 \eta \psi(t) I(t)   +(\uu/2-\eta {\bar{\beta}}_1(\T+\to)  ) |x(s)|^2     \\
&  ~~   +(\hat{e}-\eta {\bar{\beta}}_1(\T+\to) )  - \lf(  \frac{\check{\a }^2}{\uu}+\eta \iota(\T+\to) \rt)  |x(s)- x(\nu(s)-\to)|^2  \\
  &~ ~ ~~  +  {\eta} J_2(t)\Big)dt+ 2  x^T(t)  g(x(t),r(t),t)  dB(t) \Big].
  \end{split}
\end{equation}
 Due to  the increasing property of $\bar{\beta}_1(y)$ in $y>0$, one observes
\be\la{q54} \eta \bar{\beta}_1(\T+\to)\leq   (\uu/2) \wedge \hat{e}.\ee
It then follows from \eqref{q18} that
     \be\la{q55}  2 \eta \psi(t) I(t) \geq -2 \eta (\check{\a}+3\uu/4) I(t)\geq -2 \eta (\check{\a}+3\uu/4)  (\T+\to)J_2(t).\ee
     Inserting \eqref{q54} and \eqref{q55} into \eqref{q53} yields
     \begin{align} \la{q56}
d U(\bar{x}_t,\bar{r}_t, t)&\geq e^{-2  \int_0^t \psi(s) ds } \Big[ - \lf(  \frac{\check{\a }^2}{\uu}+\eta \iota(\T+\to) \rt)  |x(s)- x(\nu(s)-\to)|^2 dt \nn\\
  &  + \Big ({\eta}-2 \eta \Big(\check{\a}+\frac{3\uu}{4}\Big)  (\T+\to) \Big)J_2(t)dt+ 2  x^T(t)  g(x(t),r(t),t)  dB(t) \Big].
\end{align}
 Integrating \eqref{q56} on both sides, taking the conditional expectation with respect to $\Gt$  and using \eqref{q19}, we arrive at
\begin{align*}
 & \E_{\Gt}\lf(e^{-2  \int_0^t \psi(s) ds }\lf( \frac{\hat{e}}{2\check{d}}+ |x(t)|^2\rt) \rt)-  \eta\E_{\Gt}  \lf(e^{-2  \int_0^t \psi(s) ds }I(t) \rt)\nn\\
&\geq   \frac{ \hat{e}}{2\check{d}}+|x_0|^2-\eta  I(0)  +  \lf(  {\eta}  -\frac{2}{\uu}\check{\a }^2- \eta { {\bar{\beta}}}_3(\T+\to) \rt)   \int_0^t e^{-2  \int_0^s \psi(z) dz } \E_{\Gt} [J_2(s)]ds,
\end{align*}
where for  $y\geq 0$, \be\la{q57}  { {\bar{\beta}}}_3(y)  :=2 \iota (y) +2  y( \check{\a} +3\uu/4)=   y(8\check{\a}^2 y  +2 \check{\a}+3\uu/2 ). \ee
 Then $ \eta { {\bar{\beta}}}_3(\T+\to)\leq \uu.$  This together with the above inequality implies
\begin{align*}  \frac{\hat{e}}{2\check{d}}+ \E_{\Gt}|x(t)|^2 &\geq \lf( \frac{\hat{e}}{2\check{d}}+ |x_0|^2-\eta  I(0)\rt) e^{2  \int_0^t \psi(s) ds } .
\end{align*}
Due to \eqref{q17} and \eqref{q18} one observes from \eqref{q20} that $I(0)\leq (\T+\to) J_2(0)\leq    {\bar{\beta}}_4(\T+\to),$
where \be\la{q58}
\bar{\beta}_4(y) =y^2[3y( \bar{K}^2+\check{\a}^2)|x_0|^2 +2\bar{K}^2|x_0|^2 +\bar{K}^2 (3y+2)].\ee
Taking expectation
on both sides yields
\begin{align}\la{q59}
\frac{\hat{e}}{2\check{d}}+\E |x(t)|^2  &\geq  \lf( \frac{\hat{e}}{2\check{d}}+ |x_0|^2- \eta  {\bar{\beta}}_4(\T+\to) \rt) \E\lf(e^{2  \int_0^t  \psi(s) ds}\rt)     .
\end{align}
It follows from the definition of  $\T_*^1$ that
 $\eta { {\bar{\beta}}}_4(\T+\to)< {\hat{e}}/{2\check{d}}+|x_0|^2 $.
By virtue of Lemma \ref{l4}, for $ \ep =\uu/8 $, there is a constant $T >0$ such that
 \begin{equation*}
\E \lf( e^{2 \int_0^t  (d(r(z)))-\a(r(\nu(z)-\to))    )ds}\rt) \geq K_1(2, \a-d) e^{       7\uu t/4},~~~~ t \geq T ,
\end{equation*}
which implies
\be\la{q61}\E\lf(e^{2  \int_0^t  \psi(z)dz}\rt) \geq  K_1(2, \a-d)e^{    \uu t/4},~~~~ t \geq T .\ee
Inserting the above inequality into \eqref{q59}, we obtain
\begin{align}\la{q62}
\frac{\hat{e}}{2\check{d}}+\E |x(t)|^2  &\geq  K_1(2, \a-d)\lf( \frac{\hat{e}}{2\check{d}}+ |x_0|^2- \eta  {\bar{\beta}}_4(\T+\to) \rt) e^{   \uu t/4},~~~~ t \geq T .
\end{align}
Then the  required assertion \eqref{q43} follows.
\eproof

\subsection{Stabilization}

This subsection is to discuss  the stability and instability of DCSDS \eqref{q4} and gives the corresponding criteria.
 We replace conditions \eqref{q20} and \eqref{q21} by the following assumption in order for the SDS \eqref{q1} to
has the trivial solution $x(t)\equiv 0$.

\begin{assp} \label{A4}
 There exist positive constants  $\bar{K } $   and  $ D_i $     such that \begin{equation}\label{q63}
 |f(x, i,t)|\vee    |g(x, i,t)| \le \bar{K} |x| ,
 \end{equation} and
 \begin{equation}\label{q64}
 x^T f(x, i,t) + \frac{  1}{2} |g(x,i,t)|^2 \le D_i |x|^2
 \end{equation}
 hold for all $(x,i,t) \in \RR^n\K \Se \K \RR_+$.
 \end{assp}

 Under  Assumption \ref{A4} we will design the feedback control for the controlled system \eqref{q4} to be exponentially stable in both  mean square and almost surely (a.s.).
\begin{thm}\la{t3}
Let   Assumption \ref{A4} hold and assume that $\pi\a>\pi  D $  and  $ \ka_{\a-D}>2 $. For any $ 0<\sigma<  \zeta_{2,\a-D}^{\tilde{\T}}$, if  $0<\T\leq {\tilde{\T}} $  and $ 0<\T+\to  \leq \T_2^*(\s):= {y}_6(\s)\wedge {y}_7(\s) $,
  then
  the solution of DCSDS \eqref{q4} with  the initial solution \eqref{q7} has the properties that
\be\la{q66}
\limsup_{t\rightarrow \8}\frac{1}{t}\log\E |x(t)|^2  \leq  - (\zeta_{2,\a-D}^{ {\T}} -\s),~~
\ee and
    \be\la{q67}
 \limsup_{t\rightarrow \8}\frac{1}{t}\log( |x(t)|) \leq  -\frac{1}{2} (\zeta_{2,\a-D}^{ {\T}} -\s), ~~\PP-\hbox{a.s.} \ee
where $y_6 $ and  $y_7 $ are the positive solutions of $ {\tilde{\beta}}_1(y)=\frac{\s^2}{ 8\check{\a}^2+\s^2 }$
and
$ {\tilde{\beta}}_2(y)=\frac{\s^2}{ 8\check{\a}^2+\s^2 }$, while ${\tilde{\beta}}_1(\cdot)$ and ${\tilde{\beta}}_2(\cdot)$ are defined by \eqref{q70} and \eqref{q74}, respectively.
\end{thm}

\noindent{\it Proof.} For any $ 0<\sigma<  \zeta_{2,\a-D}^{\tilde{\T}}$, let $0<\T\leq {\tilde{\T}} $  and  $ 0<\T+\to  \leq \T_2^*(\s) $.   By the It$\hat{\hbox{o}}$ formula, the elementary inequality, and \eqref{q64}, we have
  \begin{align}
 d|x(t) |^{2}&\leq    \Big[ (2D(r(t))  - 2 \a(r(\nu(t)-\to))+\frac{\s}{2})|x(t)|^2+ \frac{2\a^2}{\s} |x(t)- x(\nu(t)-\to)|^2\Big]dt \nn\\
&~~~~~~~~~~~~+  2 x^T(t)  g(x(t),r(t),t)  dB(t).\la{q68}
 \end{align}
Using  \eqref{q63} we compute $I(t)$ and $J_1(t)$ defined by \eqref{q14} and \eqref{q16} to get
  \begin{align} \la{q69}
d I(t)
  &\leq     \tilde{\beta}_1(\T+\to)    |x(t)|^2dt    +\kappa(\T+\to)    |x(t)- x(\nu(t)-\to)|^2dt -J_2(t)dt,
\end{align}
where, for $y\geq 0$,
 \begin{align}\tilde{\beta}_1(y)   : =  y \lf[ 3y (  \bar{K}^2+ \check{\a}^2)+ \bar{K}^2\rt],~~~~~\kappa (y)  :=3 \check{\a}^2y^2.\la{q70}
   \end{align}
Define $ V(\bar{x}_t,\bar{r}_t, t)=(|x(t)|^2+\eta I(t))e^{-2  \int_0^t \phi(s) ds },$ where $ \phi(s):= D(r(s)) -\a(r(\nu(s)-\to))+\s/2,$ $\eta:=\s/2+ {4\check{\a}^2}/{\s} .$
It follows from \eqref{q68} and  \eqref{q69} that for any $t\geq 0$,
\begin{equation} \la{q71}
\begin{split}
d V(\bar{x}_t,\bar{r}_t, t)&\leq e^{-2  \int_0^t \phi(s) ds } \Big[\Big(-2 \eta \phi(t) I(t)  -(\frac{\s}{2}-\eta\tilde{\beta}_1(\T+\to)  ) |x(s)|^2    \\
&  ~~    + \lf(   \frac{2\check{\a }^2}{\s}+\eta\kappa(\T+\to) \rt)  |x(s)- x(\nu(s)-\to)|^2 \\
  &~ ~ ~~  -  {\eta} J_2(t)\Big)dt+ 2  x^T(t)  g(x(t),r(t),t)  dB(t) \Big].
  \end{split}
\end{equation}
By the definition of $\T_2^*(\s)$ one observes
$ \eta  { {\tilde{\beta}}}_1(\T+\to)\leq   \s/2.$
Integrating \eqref{q71} on both sides, taking the conditional expectation with respect to the $\sigma-$algebra $\Gt$ and using \eqref{q18} and \eqref{q19} arrives at
\begin{equation} \la{q73}
\begin{split}
 & \E_{\Gt}\lf(e^{-2  \int_0^t \phi(s) ds }|x(t)|^2\rt)+  \eta\E_{\Gt}  \lf(e^{-2  \int_0^t \phi(s) ds }I(t) \rt)\\
&\leq  |x_0|^2+\eta  I(0)  -  \lf(  {\eta}  -\frac{4 \check{\a }^2}{\s}- \eta\tilde{\beta}_2(\T+\to) \rt)   \int_0^t e^{-2  \int_0^s \phi(z) dz } \E_{\Gt} [J_2(s)]ds,
\end{split}
\end{equation}
where, for $y\geq 0$,
\be\la{q74} \tilde{\beta}_2(y)  :=2\kappa(y) +2  y \check{\a} = 2\check{\a}y(3 \check{\a} y  +1 ). \ee
By the definition of $\T_2^*(\s)$ one sees
$ \eta  {\bar{\beta}}_2(\T+\to) \leq \s/2.$
This together with  \eqref{q73}  implies that
$$ \E_{\Gt}|x(t)|^2 \leq (|x_0|^2+\eta  I(0)) e^{2  \int_0^t \phi(s) ds }.
$$
Taking expectation
on both sides, we get that
\begin{align} \la{q76}\E |x(t)|^2  &\leq  [|x_0|^2+\eta I(0)] \E\lf(e^{2  \int_0^t \phi(s) ds}\rt).
\end{align}
But, it follows from Lemma \ref{l7} that
\be\la{q77}\E\lf(e^{2  \int_0^t  \phi(s)ds}\rt) \leq  K_2(2, \a-D) e^{ -   (  \zeta_{2,\a-D}^{{\T}}-\s  )t}, ~~~~ t \geq 0.\ee Combing \eqref{q76} and \eqref{q77} yields
$ \limsup_{t\rightarrow \8}\frac{1}{t}\log\E |x(t)|^2  \leq  - (  \zeta_{2,\a-D}^{{\T}}-\s  ) ,$
which implies the required assertion \eqref{q66}.
In a similar fashion as \cite [pp. 128-130, Theorem 4.2]{mao} was proved, we can get the other required assertion \eqref{q67}.\eproof

In order to study the instability we impose the following assumption.

\begin{assp} \label{A5}
  There exist positive constants  $\bar{K }>0$,    and  $ d_i $     such that
      \eqref{q63}  and
 \begin{equation}\label{q78}
 x^T f(x, i,t) + \frac{  1}{2} |g(x,i,t)|^2 \geq d_i |x|^2
 \end{equation}
 hold for all $(x,i,t) \in \RR^n\K \Se \K \RR_+$.
 \end{assp}

\begin{thm}\la{t4}
Let Assumption \ref{A5} hold and assume $\pi\a<\pi  d$. For any  $ 0<\sigma<  \pi d-\pi \a$, if  $ 0<\T+\to  <\T_*^2(\s):=y_8(\s)\wedge y_9(\s)\wedge y_{10}(\s)$,
   then
  the solution of DCSDS \eqref{q4} with  the initial condition \eqref{q7} has the property that
\be\la{q80}
\liminf_{t\rightarrow \8}\frac{1}{t}\log\E |x(t)|^2  \geq  2( \pi d-\pi \a -\s),
\ee
where $y_i(\s)$ $(i=8,9,10)$ are the positive solutions of $ {\tilde{\beta}}_1(y)=\frac{\s^2}{ 2\check{\a}^2+\s^2 },$ $  {\tilde{\beta}}_3(y)=\frac{\s^2}{ 2\check{\a}^2+\s^2 }, $  $  {\tilde{\beta}}_4(y)=\frac{\s |x_0|^2}{ 2\check{\a}^2+\s^2 }, $ respectively, while ${\tilde{\beta}}_1(\cdot)$, ${\tilde{\beta}}_3(\cdot)$, ${\tilde{\beta}}_4(\cdot)$ are defined by \eqref{q70}, \eqref{q48}, \eqref{q49}.
\end{thm}
\noindent{\it Proof.}
For any $0<\s< \pi d -\pi \a$,  let $ 0<\T+\to  \leq \T_*^2(\s) $. Using  the elementary inequality and  \eqref{q78}, one has
  \begin{align}
 d|x(t) |^{2}
&\geq  \Big[ (2d(r(t))  - 2 \a(r(\nu(t)-\to))-\s)|x(t)|^2- \frac{\check{\a}^2}{\s} |x(t)- x(\nu(t)-\to)|^2\Big]dt \nn\\
&~~~~~~~~~~~~+  2 x^T(t)  g(x(t),r(t),t)  dB(t).\la{q81}
 \end{align}
 Define
  $\bar{U}(\bar{x}_t,\bar{r}_t, t)=( |x(t)|^2-\eta I(t))e^{-2  \int_0^t \psi(s) ds },$  where    $ \psi(s):= d(r(s)) -\a(r(\nu(s)-\to))-\s$, and $\eta:=\s+ {2\check{\a}^2}/{\s} .$
In a similar way as Theorem \ref{t2} was proved we can obtain from \eqref{q69} and \eqref{q81} that for any $t\geq 0$
\begin{equation} \la{q83}
\begin{split}
d \bar{U}(\bar{x}_t,\bar{r}_t, t)&\geq e^{-2  \int_0^t \psi(s) ds } \Big[ \Big( 2 \eta \psi(t) I(t)   +(\s-\eta {\tilde{\beta}}_1(\T+\to)  ) |x(t)|^2    \\
&  ~~    - \lf(  \frac{\check{\a }^2}{\s}+\eta \kappa(\T+\to) \rt)  |x(t)- x(\nu(t)-\to)|^2  \\
  &~ ~ ~~  +  {\eta} J_2(t)\Big)dt+ 2  x^T(t)  g(x(t),r(t),t)  dB(t)\Big].
  \end{split}
\end{equation}
One notices that
$   \eta  {\tilde{\beta}}_1(\T+\to)\leq \s$. Integrating \eqref{q83} on both sides, taking the conditional expectation with respect to the $\sigma-$algebra $\Gt$ and using \eqref{q18} and \eqref{q19}, we arrive  at
\begin{align*}
 & \E_{\Gt}\lf(e^{-2  \int_0^t \psi(s) ds } |x(t)|^2 \rt)-  \eta\E_{\Gt}  \lf(e^{-2  \int_0^t \psi(s) ds }I(t) \rt)\nn\\
&\geq  |x_0|^2-\eta  I(0)  +  \lf(  {\eta}  -\frac{2}{\s}\check{\a }^2- \eta { {\tilde{\beta}}}_3(\T+\to) \rt)   \int_0^t e^{-2  \int_0^s \psi(z) dz } \E_{\Gt} [J_2(s)]ds,
\end{align*}
where, for $y\geq 0$,
\be\la{q48}  { {\tilde{\beta}}}_3(y)  :=2 \kappa(y) +2  y( \check{\a} +\s)= 2 y(3\check{\a}^2 y  + \check{\a}+\s ).
\ee
We also see  that $\eta { {\tilde{\beta}}}_3(\T+\to)\leq   \s$. This together with the above inequality implies
\begin{align*}  \E_{\Gt}|x(t)|^2 &\geq \lf( |x_0|^2-\eta  I(0)\rt) e^{2  \int_0^t \psi(s) ds } .
\end{align*}
Due to \eqref{q17} and \eqref{q18} one observes that $I(0)\leq (\T+\to) J_2(0)\leq    {\tilde{\beta}}_4(\T+\to),$
where \be\la{q49}
\tilde{\beta}_4(y) =y^2[2y( \bar{K}^2+\check{\a}^2)|x_0|^2 + \bar{K}^2|x_0|^2  ].\ee
One notices  that $\eta { {\tilde{\beta}}}_4(\T+\to)< |x_0|^2.$
Taking expectation
on both sides yields
\begin{align}\la{q51}
\E |x(t)|^2  &\geq  \lf( |x_0|^2- \eta  {\tilde{\beta}}_4(\T+\to) \rt) \E\lf(e^{2  \int_0^t  \psi(s) ds}\rt)     .
\end{align}
 By Lemma \ref{l4}, for $0< \ep<  \pi d-\pi \a-\s $, there is a constant $T >0$ such that
\begin{align*}
\E |x(t)|^2  &\geq  \lf(   |x_0|^2-\eta   {\tilde{\beta}}_4(\T+\to) \rt)e^{    2 ( \pi d-\pi \a  -\s -\ep )t},~~~~ t \geq T .
\end{align*}
 Letting $t\rightarrow \8$, we have
 $ \liminf_{t\rightarrow \8}\frac{1}{t}\log\E |x(t)|^2  \geq  2( \pi d-\pi \a -\s-\ep).$
 As $\ep>0$ is arbitrary, the  required assertion \eqref{q80} must hold.\eproof

\section{Control of Highly Nonlinear Systems}\la{s4}
The main aim of this section is to give the easily implementable control criterion for highly nonlinear SDS \eqref{q1} such that they  stabilize \eqref{q1} exponentially in $p$th moment and almost surely. In the following, the moment and sample Lyapunov exponents are estimated, the lower bound on $\T+\to$ is given explicitly.

\subsection{Uniform Moment Boundness}
Firstly we investigate the  uniform moment boundedness of DCSDS \eqref{q4}. Generally,  SFDEs have significantly different dynamical behaviors from the corresponding SDSs.  Hence the uncontrolled  SDS \eqref{q1} may possess some property while the DCSDS \eqref{q4} may not.  We impose the following Khasminskii-type condition to guarantee that the global solution of the SDS \eqref{q1} is uniformly bounded in $p$th moment on infinite time horizon.

\begin{assp} \label{A6}
There exist positive constants $  A , B, C $ and $p\geq 2, \theta>2 $ such that
$$
 x^T f(x, i,t) + \frac{p-1}{2} |g(x,i,t)|^2 \le C+A |x|^2 -B |x|^\theta
$$
 for all $(x,i,t) \in \RR^n\K \Se \K \RR_+$.  \end{assp}

By constructing $V(x,i,t)=|x|^p$ for all $(x,i,t) \in \RR^n\K \Se \K \RR_+$, and using \cite[p. 157, Theorem 5.2]{my} we can get the following result directly. To avoid the duplication we omit the proof details.

\begin{thm} \label{t5}
Under Assumption \ref{A6}, the solution $x(t)$ of SDS (\ref{q1}) with the initial data $(x(0), r(0))=(x_0,i_0)\in \RR^n\times \Se$  satisfies
$
  \sup_{0 \le t <\8} \E|x(t)|^p < \8.
$
\end{thm}

For the DCSDS \eqref{q4} we have the following result.

\begin{thm} \label{t6}
Under Assumption \ref{A6}, the solution $x(t)$ of DCSDS (\ref{q4}) with the initial data \eqref{q7}  satisfies
$
  \sup_{ 0 \le t <\8} \E|x(t)|^p < \8.
$
\end{thm}

\noindent{\it Proof.}
 Using the It$\hat{\hbox{o}}$ formula and Assumption \ref{A6}, we  derives that, for any $t \geq 0 $,
\begin{eqnarray*}
 d(e^t |x(t) |^p)
&\leq &    e^t\Big[ pC |x( t)|^{p -2}+  (1+pA  ) |x( t)|^p - pB  |x( t)|^{p+\theta-2} \Big]dt\nn\\
&& ~~ + \check{\a}  p e^t|x( t)|^{p -1} | x(\nu(t)-\to)|dt \nn\\
&& ~~ +  pe^t|x(t)|^{p-2}   x^T(t)  g(x(t),r(t),t)  dB(t).
\end{eqnarray*}
Noting   that  for any $x, y\geq 0$
$$
 x^{p -2}\leq 1+ x^{p },~~\check{\a}  px^{p -1} y  \leq  {\frac{p-1}{p}}  (\check{\a}p)^{\frac{p}{p-1}} x^{p} +   {\frac{ 1}{p}} y^p=  {(p-1)}p^{\frac{1}{p-1}} (\check{\a})^{\frac{p}{p-1}} x^{p} +   {\frac{ 1}{p}} y^p,
$$
we have
\begin{eqnarray}\la{q85}
&&d(e^t |x(t) |^p)\nn\\
&\leq &  e^t \Big(  \bar{C}  +\frac{1}{p}  | x(\nu(t)-\to)|^p\Big) dt +  pe^t|x(t)|^{p-2}   x^T(t)  g(x(t),r(t),t) dB(t),
\end{eqnarray}
where   \be\la{q86}\bar{C}:= \sup_{x\in \RR_+} \lf\{pC+ \lf(1+pA+pC+(p-1) p^{\frac{1}{p-1}} {\check{\a}} ^{\frac{p}{p-1}}\rt) x^p  - p B  x^{p+\theta-2} \rt\}.\ee
Integrating \eqref{q85} from $0$ to $t,$
taking expectations, then dividing $e^t$ on both sides, we obtain
\begin{align*}
  \E |x(t) |^p
  &\leq \bar{C}+  |x_0|^p e^{-t}+ \frac{1}{p} \int_0^{t}  e^{s-t}  \E | x(\nu(s)-\to )|^p  ds\\
  &\leq \bar{C} +  |x_0|^p   +\frac{1}{p} \sup_{0\leq s\leq t}\lf(  \E | x(\nu(s)-\to )|^p\rt) \int_0^{t} e^{s-t}    ds\\
  &\leq \bar{C} + |x_0|^p +\frac{1}{p} \sup_{0\leq s\leq t}\lf(  \E | x(s )|^p\rt)  .
   \end{align*}
This  implies
$
  \sup_{0\leq s\leq t}\lf(  \E | x(s )|^p\rt) \leq  \frac{p(\bar{C} +|x_0|^p)}{p-1} .
$
Then the required assertion follows as  $t\rightarrow \8$.
  \eproof

\subsection{Stabilization}

In this subsection we pay attention to stabilize the nonlinear SDS \eqref{q1} by the delay feedback control based on discrete-time observations.
In order to have  the equilibrium state $0$  we further impose the following assumption.
  \begin{assp} \label{A7}
   Assume that there exist positive constants $K $, $q_1\geq 1$, $q_2\geq 1$, $ {p}\geq 2(q_1\vee q_2 ),~\theta > 2$ satisfying $ \theta  \geq  (q_1\vee q_2 )+1  $,  and $ A_i, ~ B_i $    such that
   \begin{equation}\label{q88}
 |f(x, i,t)| \le K ( |x| + |x|^{q_1}) ,~~~~
  |g(x, i,t)| \le K ( |x| +  |x|^{q_2})
 \end{equation}
    and
 \begin{equation}\label{q89}
 x^T f(x, i,t) + \frac{ p-1}{2} |g(x,i,t)|^2 \le  A_i |x|^2 -B_i  |x|^{\theta }
 \end{equation}
 hold for all $(x,i,t) \in \RR^n\K \Se \K \RR_+$.
 \end{assp}

 \begin{thm}\la{T4.4} Let Assumption \ref{A7} hold and assume that $\pi \a>\pi A$  and $ \ka_{\a-A}>2 $, where $A=(A_1,\cdots, A_N)^T$. For any $ 0<\sigma<  \zeta_{2,\a-A}^{\T'} \wedge (2 \hat{B})$ ($\T':=  \bar{\T}(2,\a-A)/2$), if  $0<\T\leq  \T'   $    and $ 0<\T+\to < \T^{**}(\s):=\bar{y}_1(\s)\wedge \bar{y}_{2}(\s) \wedge \bar{y}_3(\s) $,
 then  the solution of DCSDS \eqref{q4} with  the initial condition \eqref{q7} has the properties that
\be\la{q92}\begin{split}
&\limsup_{t\rightarrow \8}\frac{1}{t}\log\E |x(t)|^2  \leq  -  (\zeta_{2,\a-A}^{\T} -\s),\\ &\limsup_{t\rightarrow \8}\frac{1}{t}\log\E |x(t)|^\r  \leq  - (\zeta_{2,\a-A}^{\T}  -\s),
\end{split} \ee
and  \be\la{q93} \int_0^\8  \E |x(t)|^{\r +\theta-2}  dt<\8,
\ee where  $ \r:= p\wedge \theta ,$ $\bar{y}_i(\s)$ $(i=1,2,3)$ are the positive solutions of $$  2\va\beta_1(y)=\s ,~~~~  {2\va {\beta}}_2(y)=  \r(2\hat{B} -\s) , ~~~~ 2\va\beta_3(y)=\s , $$ respectively,
$\va :=\s/2+  {\check{\a}^2}[ (5\r+4)\s+ 8(\r-2)\check{A}]/ ({ \s^2}),$
${ {\beta}}_1(\cdot)$, ${ {\beta}}_2(\cdot)$, ${ {\beta}}_3(\cdot)$ are defined by \eqref{q101} and  \eqref{q114}.
 \end{thm}

\noindent{\it Proof.}  For any $ 0<\sigma< \zeta_{2,\a-A}^{\T'} \wedge (2 \hat{B})$, let $0<\T\leq  \T'$ and  $ 0<\T+\to \leq \T^{**}(\s) $.
Using the It$\hat{\hbox{o}}$ formula, the elementary inequality and \eqref{q89}, we derive
\be\la{q94}\begin{split}
 &d|x(t) |^{\r}\\
&\leq  \Big[  \r \Big(A(r(t))-\a(r(\nu(t)-\to))\Big)|x( t)|^\r    -\r (\hat{B}-\frac{\s}{2})  |x( t)|^{\r+\theta-2 }  +   \frac{\r \s}{2} |x( t)|^{2 }  \\
&  ~~+ \frac{\r \check{\a  }^2}{ 2 \s } |x(t)- x(\nu(t)-\to)|^2  \Big] dt  +    \r|x(t)|^{\r-2} x^T(t)  g(x(t),r(t),t)  dB(t).
\end{split}\ee
 In order to control the terms  $|x( t)|^\r$ and $|x( t)|^2$ together, we also derive by the It\^o formula again that
\begin{align}\la{q96}
 d|x(t) |^2
 &\leq  \dis 2\Big(A(r(t))-\a(r(\nu(t)-\to))+\frac{\s}{4} \Big) |x( t)|^2 dt  -2\hat{B}  |x( t)|^{\theta } dt\nn\\
&  ~~ + \frac{2\check{\a }^2}{\s}|x(t)- x(\nu(t)-\to)|^2   dt +   2  x^T(t)  g(x(t),r(t),t)  dB(t).
\end{align}
Under the condition \eqref{q88} we recompute $I(t)$ and $J_1(t)$ defined by \eqref{q14} and \eqref{q16}
\begin{align} \la{q97}
 J_1(t)
  &\leq (\T+\to)  \lf[  \lf( 3(\T+\to) ( K^2+ 2\check{\a}^2)+2K^2\rt) |x(t)|^2 \rt.\nn\\
  &~~\lf.+6\check{\a}^2(\T+\to) |x(t)-x(\nu(t)-\to) |^2 +3 K^2(\T+\to) |x(t)|^{2q_1}+   2K^2|x(t)|^{2q_2} \rt]\nn\\
  &\le (\T+\to)   \lf[  \lf( 6(\T+\to)(  K^2+  \check{\a}^2)+4K^2\rt) |x(t)|^2 \rt. \nn\\
  &~~~ ~\lf.  +K^2(3(\T+\to)   +   2)|x(t)|^{\r+\theta-2} +6 \check{\a}^2(\T+\to) |x(t)-x(\nu(t)-\to) |^2\rt],
\end{align}
where we have used $\r+\theta-2 \geq 2(q_1\vee q_2)$. By \eqref{q15} one has
\be \la{q100}
\begin{split}
d I(t)
  &\leq     \beta_1(\T+\to)    |x(t)|^2dt + \beta_2(\T+\to)  |x(t)|^{\r+\theta-2}dt \\  & ~~~~~~~~~~~~+6 \check{\a}^2(\T+\to)^2  |x(t)- x(\nu(t)-\to)|^2dt -J_2(t)dt,
\end{split}\ee where, for
$y \geq 0$,
 \begin{align}  \la{q101}\beta_1(y)   : = 2y\lf[ 3y (  K^2+ \check{\a}^2)+2K^2\rt], ~~
 \beta_2(y)  : =K^2y[3y   +   2] .
   \end{align}
 Define
$\bar{V}(\bar{x}_t,\bar{r}_t, t)= |x(t)|^{\r}+ \lambda |x(t)|^2+\vartheta I(t) ,$
where $ \lambda:= 1+\r+2(\r-2)\check{A}/\s $, and $\va$ is given in the theorem.  For any $t\geq 0$, define
$ \varphi(t):=A(r(t))-\a(r(\nu(t)-\to))+\s/2.$
Using \eqref{q94}, \eqref{q96} and \eqref{q100} arrives at
\begin{align} \la{q103}
&d \bar{V}(\bar{x}_t,\bar{r}_t, t)  \nn  \\
 &\leq  \Big[ \r (A(r(t))-\a(r(\nu(t)-\to)))|x( t)|^\r    -(\r \hat{B}- \frac{\r\s}{2}-\va  \beta_2(\T+\to) )  |x( t)|^{\r+\theta-2 }  \nn  \\
&  ~~~~~~~~~~~~ - 2\lambda \hat{B}  |x( t)|^{\theta } + \dis  \lf( 2\lambda \varphi(t)+  \frac{\r\s}{2}- \frac{\lambda\s}{2}  +\va  \beta_1(\T+\to) \rt) |x( t)|^2   -\va J_2(t)\nn \\
&   ~~~~~~~~~~~~     + \Big(\frac{ (\r+4\lambda) \check{\a  }^2}{ 2 \s }  +6\va  \check{\a}^2(\T+\to)^2  \Big) |x(t)- x(\nu(t)-\to)|^2  \Big] dt \nn \\
&   ~~~~~~~~~~~~ +    (\r |x(t)|^{\r-2}+2\lambda) x^T(t)  g(x(t),r(t),t)  dB(t).
\end{align}
 Since $\varphi(s)$ has a finite number of jumps in any  finite interval $[0, t]$, $\int_0^t  \varphi(s)ds$ is derivable.
 Thus, it follows from \eqref{q103} that
\be \la{q104}
\begin{split}
&d \lf(e^{-2\int_0^t  \varphi(s)ds }\bar{V}(\bar{x}_t,\bar{r}_t, t)  \rt)   \\
 &=  \dis e^{-2 \int_0^t  \varphi(s) ds }\Big(-2 \varphi(t)\bar{V}(\bar{x}_t,\bar{r}_t, t)dt+d \bar{V}(\bar{x}_t,\bar{r}_t, t)   \Big) \\
 &=  e^{-2 \int_0^t \varphi (s) ds }\Big[ - 2  \va \varphi (t)I(t)+(\r-2) (A(r(t))-\a(r(\nu(t)-\to)))|x( t)|^\r \\
&~~~~~~~~~~~~~~~~ ~~~- \s|x( t)|^\r-(\r \hat{B}-\frac{\r\s}{2}-\va  \beta_2(\T+\to) )  |x( t)|^{\r+\theta-2 }  \\
&~~~~~~~~~~~~~~~~ ~~~- 2\lambda \hat{B}  |x( t)|^{\theta }  - \dis  \lf( \frac{(\lambda-\r)\s}{2} -\va  \beta_1(\T+\to) \rt) |x( t)|^2  -\va J_2(t)  \\
&~~~~~~~~~~~~~~~~ ~~~     + \Big(\frac{ (\r+4\lambda) \check{\a  }^2}{ 2\s }  +6\va  \check{\a}^2(\T+\to)^2  \Big) |x(t)- x(\nu(t)-\to)|^2  \Big] dt \\
           &~~~~~~~~~~~~~~~~ ~~~+      e^{-2 \int_0^t \varphi (s) ds }(\r |x(t)|^{\r-2}+2\lambda) x^T(t)  g(x(t),r(t),t)  dB(t) .\end{split}\ee
  One observes from \eqref{q18} that
     \be\la{q105} -2 \va \varphi(t) I(t) \leq 2 \va \check{\a} I(t)\leq 2 \va\check{\a}  (\T+\to)J_2(t).\ee
    Noticing $2\leq \r\leq \theta $,  we obtain
\be\la{q106}
 (A(r(t))-\a(r(\nu(t)-\to)))|x( t)|^\r\leq \check{A}|x( t)|^\r\le \check{A}|x(t)|^2+\check{A} |x(t)|^{ \theta } . \ee
  Inserting \eqref{q105}  and \eqref{q106} into \eqref{q104} yields
   \begin{align}\la{q107}
     &d \lf(e^{-2\int_0^t  \varphi(s)ds }\bar{V}(\bar{x}_t,\bar{r}_t, t)  \rt)  \nn\\
           &\leq     e^{-2 \int_0^t \varphi (s) ds }\Big[  -(\r \hat{B}-\frac{\r \s}{2}-\va  \beta_2(\T+\to) )  |x( t)|^{\r+\theta-2 }  - (2\lambda \hat{B}-(\r-2)  \check{A})  |x( t)|^{\theta }\nn\\
  &~~~~~~~~~~  - \dis  \lf( \frac{(\lambda -\r) \s}{2}-(\r-2)  \check{A} -\va  \beta_1(\T+\to) \rt) |x( t)|^2   -\va(1-2  \check{\a}  (\T+\to) )  J_2(t)  \nn\\
  &~~~~~~~~~~   + \Big(\frac{ (\r+4\lambda) \check{\a  }^2}{ 2 \s }  +6\va  \check{\a}^2(\T+\to)^2  \Big) |x(t)- x(\nu(t)-\to)|^2  \Big] dt\nn\\
           &~~~~~~~~~~~~~~~~ ~~~+      e^{-2 \int_0^t \varphi (s) ds }(\r |x(t)|^{\r-2}+2\lambda) x^T(t)  g(x(t),r(t),t)  dB(t).
\end{align}
 One notices from $0<\sigma< 2 \hat{B}$ and the definition of $\lambda $ that
  \be\la{q108}
 2\lambda \hat{B}-(\r-2)  \check{A}>0.
 \ee
  From the definitions of $\T^{**}(\s)$, $\lambda$, $\beta_1(\cdot)$ and $\beta_2(\cdot)$, one notices that $\T^{**}(\s)<1$, and furthermore
\begin{align}
  \va\beta_1(\T+\to) \leq \frac{\s}{2}=  \frac{\s(\lambda  - \r)}{ 2}-(\r-2)  \check{A}, ~~~~
   \va\beta_2(\T+\to) \leq \r( \hat{B}-  \frac{\s}{2}).\la{q110}
\end{align}
 Substituting \eqref{q108}-\eqref{q110} into \eqref{q107} yields
 \begin{equation}\la{q111}
 \begin{split}
     &d \lf(e^{-2\int_0^t  \varphi(s)ds }\bar{V}(\bar{x}_t,\bar{r}_t, t)  \rt)  \\
           &\leq     e^{-2 \int_0^t \varphi (s) ds }\Big[   \Big(\frac{ (\r+4\lambda) \check{\a  }^2}{ 2\s }  +6\va  \check{\a}^2(\T+\to)^2  \Big) |x(t)- x(\nu(t)-\to)|^2 dt  \\
  &~~~~   -\va(1-2  \check{\a}  (\T+\to) )  J_2(t)   dt +     (\r |x(t)|^{\r-2}+2\lambda) x^T(t)  g(x(t),r(t),t)  dB(t) \Big].
  \end{split}
\end{equation}
  Using \eqref{q18} implies that \be\la{q112} \bar{V}(\bar{x}(0),\bar{r}(0), 0)\leq |x_0|^{\r}+ \lambda |x_0|^2+\vartheta (\T+\to)J_2(0)<\8. \ee
Integrating \eqref{q111} on both sides, taking the conditional expectation with respect to the $\sigma-$algebra $\Gt$ and using \eqref{q112},  \eqref{q19}, we  obtain
\begin{align} \la{q113}
& e^{-2 \int_0^t  \varphi(s)ds }\E_{\Gt}\bar{V}(\bar{x}_t,\bar{r}_t, t) \leq  \bar{V}(\bar{x}(0),\bar{r}(0), 0)  \nn\\
  &~~~ -\int_0^t e^{-2 \int_0^s \varphi (z) dz }\Big[
              (\va- \frac{ (\r+4\lambda) \check{\a  }^2}{ \s } -\va \beta_3(\T+\to)  )  \E_{\Gt}[J_2(s)]\Big]ds ,
\end{align}
where for any $y\geq 0$,
\be\la{q114} \beta_3(y) :=2\check{\a}y(1+6 \check{\a} y).\ee
Due to $p\geq 2(q_1 \vee q_2)$ one observes  from Theorem \ref{t6} that $\E_{\Gt}[J_2(s)]<\8$ for any $s\geq 0$.
It follows from the definitions of $\T^{**}(\s)$, $\lambda$,  $\va$  that
$
   \va \beta_3(\T+\to) \leq  \frac{\s}{2} = \va- \frac{ (\r+4\lambda) \check{\a  }^2}{  \s }.
$
This together with \eqref{q113} implies
\begin{align*}& \E_{\Gt} \bar{V}(\bar{x}_t,\bar{r}_t, t)   \leq   \bar{V}(\bar{x}(0),\bar{r}(0), 0)e^{ 2 \int_0^t  \varphi(s)ds }  ds.\end{align*}
Then
$$  \E |x(t)|^{\r}+ \lambda \E |x(t)|^2 \leq \bar{V}(\bar{x}(0),\bar{r}(0), 0)  \E\lf( e^{ 2 \int_0^t  \varphi(s)ds }\rt).$$
 It follows from Lemma \ref{l7} that
\be\la{q115}\E |x(t)|^{\r}+ \lambda\E |x(t)|^2 \leq  \bar{V}(\bar{x}(0),\bar{r}(0), 0) e^{ -  (\zeta_{2,\a-A}^{\T} -\s)t },\ee  which implies that \eqref{q92} holds.
 Integrating  \eqref{q103} on both sides, taking expectation, and using the similar techniques yields
 \begin{align*}
   & (\r \hat{B}-\frac{\r \s}{2}-\va  \beta_2(\T+\to) )   \int_0^t  \E  |x( s)|^{\r+\theta-2 } ds \\
   &\leq  \bar{V}(\bar{x}(0),\bar{r}(0), 0) + \r(\check{A}+\frac{  \s}{2}) \int_0^t   ( \E |x( s)|^\r + \lambda \E |x( s)|^2)ds .
\end{align*}
This, together with \eqref{q115},  implies
 $$
     (\r \hat{B}-\frac{\r \s}{2}-\va  \beta_2(\T+\to) )   \int_0^t  \E  |x( s)|^{\r+\theta-2 } ds
    \leq  \bar{C}_1 ,$$ where $\bar{C}_1$ is a positive constant. The conclusion \eqref{q93} follows by
 letting $t\rightarrow \8.$
 \eproof

The corresponding results for a  special case $p\geq \theta$ follows directly from the above proof but holds for a possible bigger $\T^{**}(\s)$.

\begin{cor}\la{c2} Let Assumption \ref{A7} hold with $p\geq \theta$, $\pi \a>\pi A$  and $ \ka_{\a-A}>2 $. For any $ 0<\sigma<  \zeta_{2,\a-A}^{\T'} \wedge (2 \hat{B})$, if  $0<\T\leq  \T' $ and  $ 0<\T+\to < \bar{\T}^{**} (\s):= \bar{y}_1'(\s)\wedge\bar{y}_2'(\s)\wedge\bar{y}_3'(\s)$, the conclusions of Theorem \ref{T4.4} hold with $\r=\theta$, where $\bar{y}_i'(\s)$ $(i=1,2,3)$ are the positive solutions of $  4\hat{B}\va_1\beta_1(y)=\sigma(\theta-2)\check{A} ,$ $ 2\va_1 { {\beta}}_2(y)=  \theta(2\hat{B} -\s) $,  $ 2\va_1\beta_3(y)=\s $, respectively,
$  \va_1 :=\s/2+ {\check{\a}^2}[   \theta \hat{B}+ 2(\theta-2)\check{A}]/{ (\s \hat{B})},$
${ {\beta}}_1(\cdot)$, ${ {\beta}}_2(\cdot)$, ${ {\beta}}_3(\cdot)$ are defined by \eqref{q101} and  \eqref{q114}.
\end{cor}

%By virtue of Theorem \ref{T4.4} we get the asymptotic stability of   $\E |x(t)|^q$ for $q\in [2, p]$ by the Cauchy inequality  if $p>\theta$.
Due to the uniform boundedness of $\E |x(t)|^p$ on the infinite horizon, by the H\"older inequality, we go one step further to obtain the following result.
\begin{thm}\la{t8}
Under the conditions of Theorem  \ref{T4.4},   for any $ 0<\sigma<  \zeta_{2,\a-A}^{\T'} \wedge (2 \hat{B})$, if  $0<\T\leq  \T' $ and  $ 0<\T+\to <\bar{\T}^{**} (\s)$,
  the solution of DCSDS \eqref{q4} with  the initial condition \eqref{q7} has the property that for any $q\in [2, p)$
  \be\la{q123}
 \limsup_{t\rightarrow \8}\frac{1}{t}\log(\E |x(t)|^q) \leq \xi_q:= \lf\{ \begin{array}{lcl}
 -  (\zeta_{2,\a-A}^{\T}-\s),~~~~~~   q =2,&\\
 -  \frac{ q }{  \r}(\zeta_{2,\a-A}^{\T}-\s),~~~ ~  q \in (2,  \r],& (\hbox{if} ~\r=p),\\
  -\frac{ p-q }{ p-\r} (\zeta_{2,\a-A}^{\T}-\s),~~  q \in(\r, p), &~(\hbox{if} ~\r<p).
  \end{array}\rt. \ee
  \end{thm}

Using the techniques of \cite[p.10, Theorem 4.5]{lm} we can obtain the following sample Lyapunov exponent. But to avoid duplication we omit the proof.

 \begin{thm}\la{t9}
 Under the conditions of Theorem  \ref{T4.4} and $p>   v:= (2q_1)\vee( 2q_2) $,   for any $ 0<\sigma<  \zeta_{2,\a-A}^{\T'} \wedge (2 \hat{B})$, if  $0<\T\leq  \T' $ and  $ 0<\T+\to \leq \T^{**} (\s)$,
  the solution  of DCSDS \eqref{q4} with  the initial solution \eqref{q7} has the property that
  \be\la{q124}
 \limsup_{t\rightarrow \8}\frac{1}{t}\log( |x(t)|) \leq  -\xi_v ~~  ~\hbox{a.s.} \ee
where the definition of $\xi_{\cdot}$ is given by \eqref{q123}. This implies that the  DCSDS \eqref{q4} is almost surely exponentially stable
 \end{thm}

\section{Example}\la{s5}
\begin{expl} Consider a scalar nonlinear SDS \eqref{q1} with  a scalar Brownian motion
 $B(t) $, a Markov chain $r(t)$ on the state space $\Se=\{1,~2\}$ with its generator matrix
$\G=\lf(\begin{array}{lcl}
 -10 &10\\
 20&-20
\end{array}\rt)$, and  the coefficients $f$ and $g$ defined by
\begin{align*}
f(x, 1,t) =x (1-3x^2 ),&~~~~g(x, 1,t) =|x |^{3/2},\\
f(x, 2,t) =x (1- 2x^2 ),&~~~~g(x, 2,t) =x.
\end{align*}\end{expl}
One observes that \eqref{q88} is satisfied with $q_1=3,~q_2=3/2, K=3$.
Due to the  Young inequality one goes a further step to obtain that \be\la{q129}\begin{split}
 x^T f(x, 1,t)+\frac{7-1}{2}|g(x, 1,t)|^2 &=|x|^2+3|x|^3-3|x|^4\leq 2.5|x|^2-1.5|x|^4  , \\
 x^T f(x, 2,t)+\frac{7-1}{2}|g(x, 2,t)|^2 &=4|x|^2 -2|x|^4,
\end{split}\ee
which implies that \eqref{q89} is satisfied with $p=7,~\theta=4, \r=4$, $ A_1=2.5, B_1=1.5, A_2=4, B_2=2$. Thus Assumption \ref{A7} holds. By a direct computation we know the stationary distribution $(\pi_1,\pi_2)=(2/3, 1/3)$ and
$ \pi A= 3$.
%Set $u(x,i,t)=-\a_i x$ for any $(x,i,t)\in \RR\times \Se\times \RR_+$.
By virtue of Theorem \ref{t6} the controlled system \eqref{q4}
%\begin{equation} \label{q130}
% dx(t)= [f(x(t),r(t),t)-\a(r(\nu(t)-\to))x(r(\nu(t)-\to)) ]dt + g(x(t), r(t),t)dB(t)
%\end{equation}
with any initial value condition \begin{align}\label{q131}
   x(t)  =x_0\in \RR, ~~ ~ ~   r(t)=i_0\in \Se,~  - \to\le t\le 0,
\end{align}
has a  unique global solution $x(t)$ on $[0,\8)$ which satisfies
$
  \sup_{0 \le t <\8} \E|x(t)|^7 < \8.
$
\begin{figure}[h!tb]
 \centering
\includegraphics[height=3cm, width=10cm]{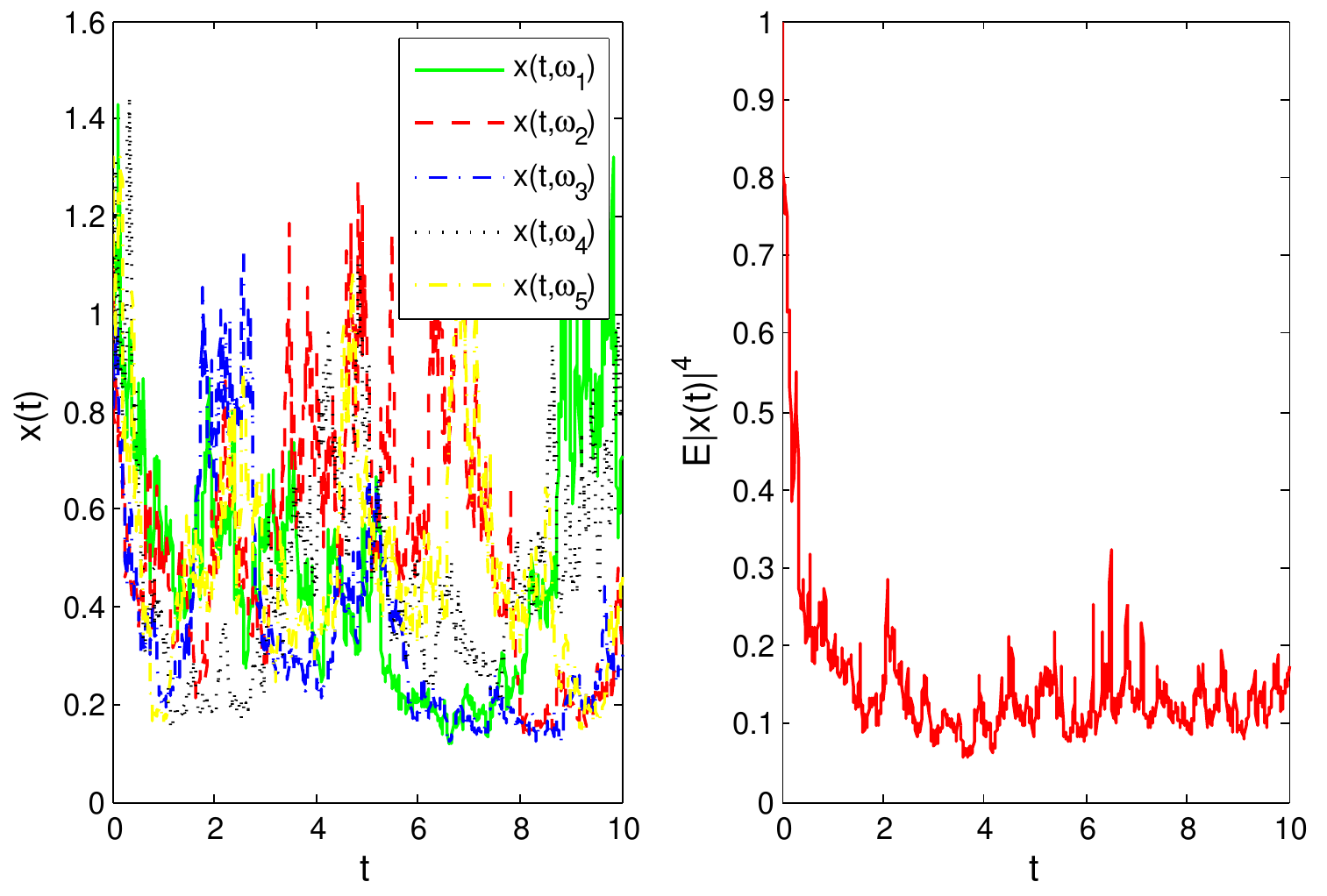}
 \caption{ Five  sample pathes of the solution $x(t)$  of \eqref{q1} and the sample mean of $|x(t)|^4$ for $100$ sample points on $t\in [0,10]$ with the initial value $(x(0),r(0))=(1,2)$ and step size $\triangle=10^{-6}$.}\label{f2}
\end{figure}

In order to have a feeling on the asymptotic behavior we carry out some numerical simulations using  MATLAB with the time step size   $\triangle=10^{-6}$. %Figure \ref{f1} depicts a sample path of the Markov chain $r(t)$ with $r(0)=2$ on $[0,10]$;
Figure \ref{f2} depicts $5$ sample pathes of the solution and the sample mean of $ |x(t)|^4$ for $100$ sample points, with the initial value $(x(0),r(0))=(1,2)$ for  $t\in [0,10]$. One observes from Figure \ref{f2} that the solutions is uniformly bounded in the $4$th moment, but the trivial solution $x(t)\equiv 0$ is  unstable  either $\PP$-a.s. or in the moment. So it is necessary to input the feedback control to stabilize  SDS \eqref{q1}. We will discuss two cases on the design of control functions. In both cases, we will give the range for $\T+\to$ to take and estimate the corresponding Lyapunov exponents.

\begin{figure}[h!tb]
 \centering
\includegraphics[height=3cm, width=10cm]{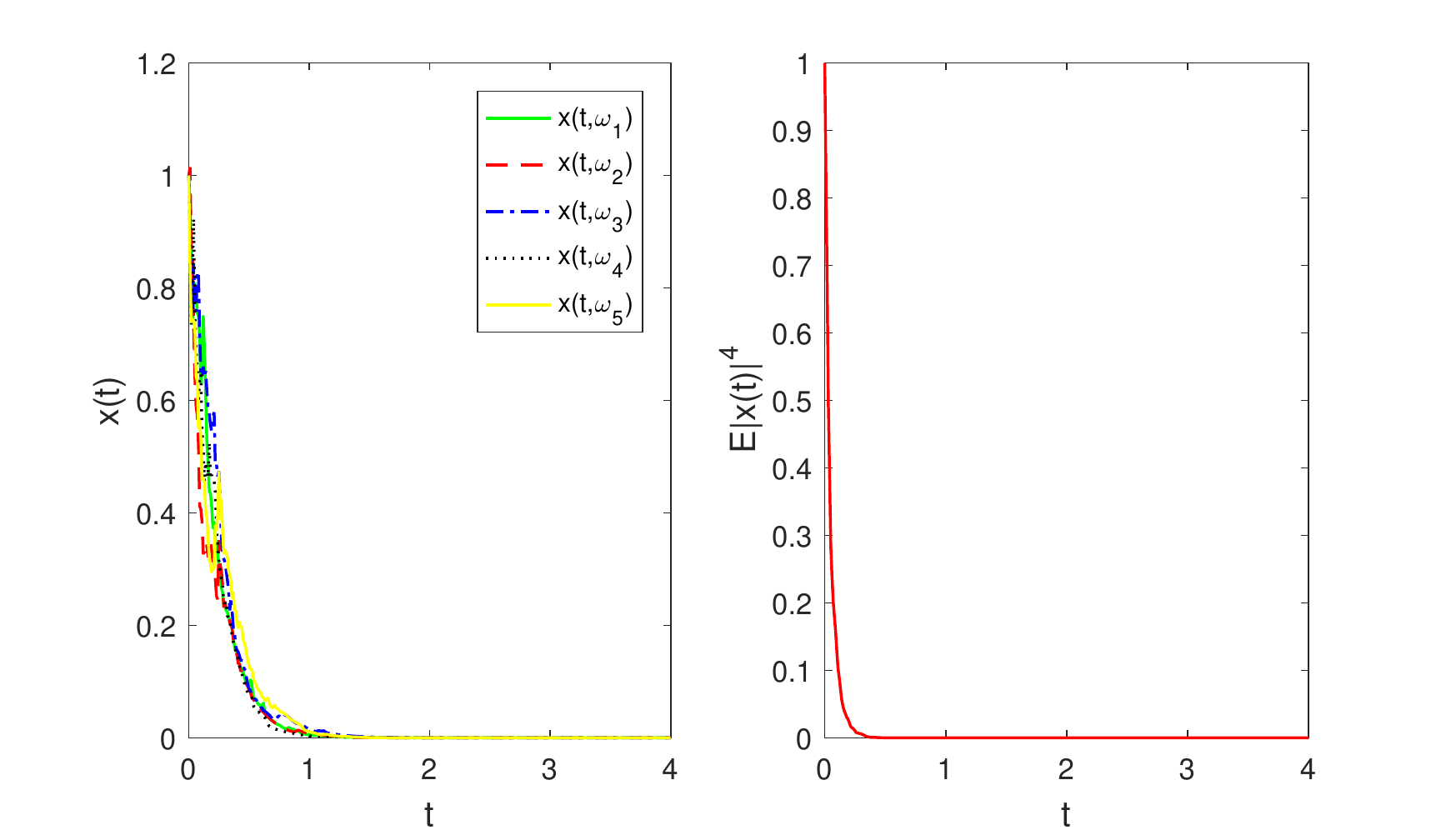}
  \caption{ Five  sample pathes of the solution $x(t)$ and the sample mean of $|x(t)|^4$ for $100$ sample points,  for the controlled system \eqref{q4} for $t\in [0,4]$ with the initial value $(x(0),r(0))=(1,2)$ and step size $\triangle=10^{-6}$.}\label{f3}
\end{figure}

\underline{Case $1$} In this case we consider that the state of the underlying SDS and
the Markov chain are  observable and the feedback control can be input in both modes $1$ and $2$.  Let $\a( 1 )= 6 $, $\a (2)= 6 $. Then $\pi\a=6>\pi A$ and $ \ka_{\a-A}=\8 $. By \eqref{q-18} and \eqref{q-18+} we can obtain that
$ {\T'}=9.6\times10^{-3}$ and $\zeta_{2,\a-A}^{\T'}=3.265 $.  Fix $\s=2<3= \zeta_{2,\a-A}^{\T'} \wedge (2\hat{B} )$, we may get $\T^{**}_1(\s)=2.78\K 10^{-4}, $ choose $\T=1\K 10^{-4},~\to= 1.7 \K 10^{-4}$, then $\zeta_{2,\a-A}^{\T}=5.8345 $. By virtue of Theorem \ref{T4.4}, the solution of DCSDS \eqref{q4} with  the initial condition \eqref{q131} has the properties that
\begin{align*}
\limsup_{t\rightarrow \8}\frac{1}{t}\log\E |x(t)|^2  \leq  -3.8345,~&~ \limsup_{t\rightarrow \8}\frac{1}{t}\log\E |x(t)|^4  \leq  -3.8345,\\
 \int_0^\8  \E |x(t)|^{6}  dt<\8,~&~\limsup_{t\rightarrow \8}\frac{1}{t}\log( |x(t)|) \leq  -1.9172  ~~\PP-\hbox{a.s.}
\end{align*}
Figure \ref{f3} depicts five sample pathes of the solution $x(t)$ and the sample mean of $|x(t)|^4$ for $100$ sample points for the controlled system \eqref{q4} for $t\in [0,4]$ with the initial value $(x(0),r(0))=(1,2)$ and step size $\triangle=10^{-6}$.

\begin{figure}[h!tb]
 \centering
\includegraphics[height=3cm, width=10cm]{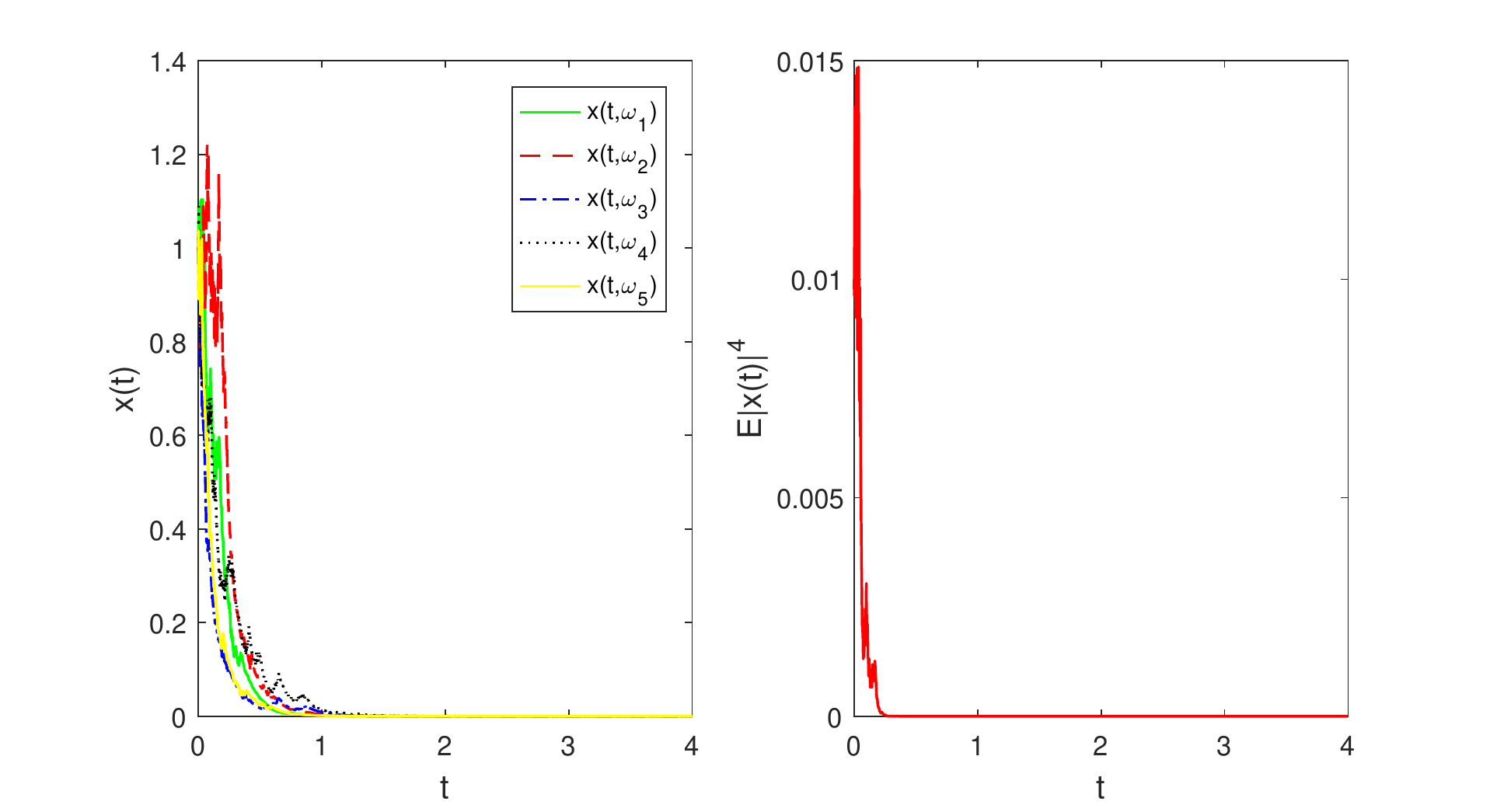}
 \caption{ Five  sample pathes of the solution $x(t)$ and the sample mean of $|x(t)|^4$  for $100$ sample points for the controlled system \eqref{q4}
  where the control is only input to the system in mode 1, for $t\in [0,4]$ with the initial value $(x(0),r(0))=(1,2)$ and step size $\triangle=10^{-8}$.}\label{f4}
\end{figure}

\underline{Case $2$} In this case we consider that the feedback control can only be input to one mode but not the other.  Assume that the system in mode $1$ is controllable but not in mode $2$. Mathematically, we let $\a( 1 )= 9  $, $\a (2)=0 $.  Then $\pi\a=6$ and $ \ka_{\a-A}=3.46 $. By \eqref{q-18} and \eqref{q-18+} we can obtain that
$ {\T'}=3.73 \times10^{-3}$ and $\zeta_{2,\a-A}^{\T'}=0.5626 $.  Fix $\s=0.5<0.5626= \zeta_{2,\a-A}^{\T'} \wedge (2\hat{B} )$, we may get $\T^{**}_2(\s)=5.83\K 10^{-6}, $ choose $\T=3\K 10^{-6},~\to= 2.8 \K 10^{-6}$, then $\zeta_{2,\a-A}^{\T}=1.0747$.
By Theorem \ref{T4.4}, we can then conclude that
\begin{align*}
\limsup_{t\rightarrow \8}\frac{1}{t}\log\E |x(t)|^2  \leq  -0.5747,~&~ \limsup_{t\rightarrow \8}\frac{1}{t}\log\E |x(t)|^4  \leq  -0.5747,\\
 \int_0^\8  \E |x(t)|^{6}  dt<\8,~&~\limsup_{t\rightarrow \8}\frac{1}{t}\log( |x(t)|) \leq  -0.2874~ ~\PP-\hbox{a.s.}
\end{align*}
Figure \ref{f4} depicts five  sample pathes of the solution $x(t)$ and the sample mean of $|x(t)|^4$  for $100$ sample points for the controlled system \eqref{q4} for
   $t\in [0,4] $ with the initial value $(x(0),r(0))=(1,2)$, step size $\triangle=10^{-8}$.
Due to the definition of $\beta_3(\cdot) $ in \eqref{q114} one observes that the balanced control values $\a(\cdot)$ in modes are helpful to get a better lower bound of $\T^*$ or $\T^{**}$, and this is illustrated in this example.

\section*{Acknowledgement}

The authors would like to thank the editor and reviewers for their very helpful comments and suggestions.


\begin{thebibliography}{17}

\bibitem{An} W.J. Anderson, Continuous Markov chains, Springer, New York, 1991.


\bibitem{ABG}  A. Arapostathis, V.S. Borkar and M.K. Ghosh,
 Ergodic control of diffusion processes, Cambridge University Press, Cambridge, 2012.

\bibitem{Ba} J.B. Bardet,  H. Gu$\acute{e}$rin, F. Malrieu,  Long time behavior of diffusions with Markov switching, ALEA Lat. Am. J. Probab. Math. Stat., 7 (2010), pp. 151-170.

    \bibitem{ds1} A.B.  Chammas and C.T. Leondes, On the finite time control of linear systems by
piecewise constant output feedback, Internat. J. Control, 30 (1979),  pp. 227-234.

%\bibitem{chen} M.F. Chen, From Markov Chains to Non-Equilibrium Particle Systems, 2nd ed., %World Scientific, River Edge, NJ, 2004.

\bibitem{ds2} T. Hagiwara  and M. Araki, Design of stable state feedback controller based on the
multirate sampling of the plant output, IEEE Trans. Automat. Control, 33 (1988), pp.
812-819.

\bibitem{ds3}  T. Hagiwara  and M. Araki, On preservation of strong stabilizability under sampling.
IEEE Trans. Automat. Control, 33 (1988), pp. 1080-1082.


\bibitem{gg} J.C. Geromel and G.W. Gabriel, Optimal H2 state feedback sampled-data control design of Markov jump linear systems, Automatica, 54 (2015), pp. 182-188.

\bibitem{lm} X. Li, X. Mao, Stabilisation of highly nonlinear hybrid stochastic differential delay equations by delay feedback control, Automatica, 112 (2020), 108657.

\bibitem{Li17}
Y. Li, J. Lu, X. Mao and Q. Qiu,
 Stabilization of hybrid systems by feedback control based on discrete-time state and mode observations,
 Asian J. Control 19(6) (2017), pp.1943--1953.

\bibitem{lu}Z. Lu, J. Hu,  and X.  Mao,  Stabilisation by delay feedback control for highly nonlinear hybrid stochastic differential equations. Discrete and Continuous Dynamical Systems - Series B.,  24 (2019), pp.  4099-4116.

    \bibitem{mao}   X. Mao,   Stochastic Differential Equations and Applications, 2nd Edition, Chichester: Horwood, 2008.

  \bibitem{mao13} X. Mao, Stabilization of continuous-time hybrid stochastic differential equations by discrete time feedback control, Automatica, 49 (2013), pp. 3677-3681.

      \bibitem{mao08} X. Mao, L. James, L. Huang, Stabilisation of hybrid stochastic differential equations by delay feedback control, System and Control Letter, 57 (2008),  pp. 927-935.

 \bibitem{mao14} X. Mao, W. Liu, L. Hu, Q. Luo, and J. Lu, Stabilization of hybrid stochastic differential equations by feedback control based on discrete-time state observations, Systems Control Lett., 73 (2014), pp. 88-95.

     \bibitem{my} X. Mao, C. Yuan, Stochastic Differential Equations with Markovian Switching,
  Imperial College Press, 2006.





\bibitem{Pin1} M. Pinsky, R.G. Pinsky, Transience/recurrence and central limit theorem
  behavior for diffusions in random temporal environments, Ann. Probab. 21
  (1993),  pp.  433--452.

\bibitem{Pin2}
M. Pinsky, M. Scheutzow, Some remarks and examples concerning the transience
  and recurrence of random diffusions, Ann. Inst. H. Poincar$\acute{\hbox{e}}$
  Probab. Statist. 28 (1992),  pp. 519--536.

  \bibitem{qiu} Q. Qiu, W. Liu, L. Hu, X. Mao, S. You, Stabilization of stochastic differential equations with Markovian switching by feedback control based on discrete-time state observation with a time delay, Statistics and Probability Letters, 112 (2016), pp. 16-26.


  \bibitem{rm} M.L. Rosinberg, T. Munakata, and G. Tarjus, Stochastic thermodynamics of Langevin systems under time-delayed feedback control: Second-law-like inequalities,
Phys. Rev. E 91 (2015), 042114.


 \bibitem{shao0} J. Shao, Strong solutions and strong Feller properties for regime-switching diffusion processes in an infinite state space, SIAM J. Control Optim., 53 (2015), pp. 2462-2479.

  \bibitem{shao} J. Shao, Stabilization of regime-switching processes by feedback control based on discrete time observations. SIAM J. Control Optim. 55 (2017),  pp. 724-740.

\bibitem{shao2}J. Shao and F. Xi, Stabilization of regime-switching processes by feedback control based on discrete time observations II: State-dependent case. SIAM J. Control Optim. 57 (2019), pp. 1413-1439.

\bibitem{Shiryaev}{A.N. Shiryaev,  {Probability  2nd edition}, World Book Inc., 2004.}


\bibitem{song} G. Song, B. Zheng, Q. Luo, X. Mao, Stabilisation of hybrid stochastic differential equations by feedback control based on discrete-time observations of state and mode, IET Control Theory Appl., 11 (2017), pp. 301-307.

\bibitem{yf10} G. Yin and F. Xi, Stability of regime-switching jump diffusions, SIAM J. Control Optim., 48 (2010), pp. 4525-4549.

\bibitem{yz09} G. Yin and C. Zhu,{Hybrid Switching Diffusions. Properties and Applications},
Springer, New York, 2010.



\bibitem{mao15} S. You, W. Liu, J. Lu, X. Mao, and Q. Qiu, Stabilization of hybrid systems by feedback control based on discrete-time state observations, SIAM J. Control Optim., 53 (2015), pp. 905-925.











\end{thebibliography}
\end{document}